\documentclass[reqno, 11pt]{amsart}

\usepackage{bbm}
\usepackage{euscript}
\usepackage{pb-diagram}
\usepackage{lamsarrow}
\usepackage{pb-lams}
\usepackage{amsmath}
\usepackage{amsthm}
\usepackage{epsfig}
\usepackage{amssymb}
\usepackage{pifont}
\usepackage{amsbsy}

\newtheorem{theorem}{Theorem}[section]
\newtheorem{lemma}[theorem]{Lemma}
\newtheorem{proposition}[theorem]{Proposition}

\theoremstyle{definition}
\newtheorem{definition}[theorem]{Definition}
\newtheorem{ex}[theorem]{Example}

\theoremstyle{remark}
\newtheorem{remark}[theorem]{Remark}

\numberwithin{equation}{section}

\newskip\aline \newskip\halfaline
\def\skipaline{\vskip\aline}

\def\qedbox{$\rlap{$\sqcap$}\sqcup$}
\def\qed{\nobreak\hfill\penalty250 \hbox{}\nobreak\hfill\qedbox\skipaline}

\def\proofend{\eqno{\mbox{\qedbox}}}

\newcommand\bC{{\mathbb C}}

\newcommand\bR{{\mathbb R}}

\newcommand\bZ{{\mathbb Z}}

\newcommand{\C}{\EuScript{C}}

\newcommand{\e}{\EuScript{E}}
\newcommand{\f}{\EuScript{F}}

\newcommand{\h}{\EuScript H}

\newcommand{\m}{\EuScript{M}}
\newcommand{\n}{\EuScript{N}}
\newcommand{\p}{\EuScript{P}}

\newcommand{\s}{\EuScript{S}}

\newcommand{\cv}{\EuScript{V}}

\newcommand{\ra}{\rightarrow}

\newcommand{\Lra}{{\longrightarrow}}

\def\inpr{\mathbin{\hbox to 6pt{\vrule height0.4pt width5pt depth0pt \kern-.4pt \vrule height6pt width0.4pt depth0pt\hss}}}

\newcommand{\ii}{{\bf i}}

\newcommand{\ve}{{\varepsilon}}
\newcommand{\vfi}{{\varphi}}

\newcommand{\pa}{\partial}

\def\mmod{\mathop{\rm \,mod}\nolimits}

\newcommand{\one}{\mathbbm{1}}

\DeclareMathOperator{\Cr}{\mathbf{Cr}}

\begin{document}

\title{ Morse functions of the two sphere}
\date{Started October, 2005. Finished  January 2006.}

\author{Liviu I. Nicolaescu}

\address{Department of Mathematics, University of Notre Dame, Notre Dame, IN 46556-4618.}
\email{nicolaescu.1@nd.edu}
\urladdr{http://www.nd.edu/~lnicolae/}

\begin{abstract}
We count how many "different" Morse functions exist on the $2$-sphere. There are several ways of declaring that two Morse functions $f$ and $g$ are "indistinguishable" but  we concentrate   only two natural equivalence relations: homological (when the regular sublevel sets  $f$ and $g$ have identical Betti numbers), and geometric (when $f$ is obtained from $g$ via global, orientation preserving  changes of coordinates on $S^2$ and $\bR$). The   count of homological classes is reduced  to  a count of lattice paths   confined to  the first quadrant. The count   of geometric classes is reduced   to a count of certain labelled trees. We produce  a two-parameter recurrence  which  can be encoded  by a first order  quasilinear  pde. We solve this  equation  using the classical method of characteristics and we produce  a closed form description of the exponential generating  function of the numbers of geometric classes.
\end{abstract}

\maketitle

\tableofcontents

\newpage

\section{The main problem}
\setcounter{equation}{0}

Suppose $X$ is a smooth compact, oriented  manifold without boundary. Following R. Thom, we say that $f: X\ra \bR$ is an \emph{excellent} Morse function if all its critical points are nondegenerate and no two of them correspond to the same critical value. We denote by $\m_X$ the space of \emph{excellent} Morse functions on $X$.  To ease the presentation,   in the sequel a Morse function will by default   be excellent.

For $f\in\m_X$ we denote by $\nu(f)$ the number  of critical points of $f$. Given a   Morse function $f: X\ra \bR$ with $\nu(f)=n$  we define a \emph{slicing} of $f$ to be an increasing sequence of real numbers
\[
-\infty=a_0<a_1<\cdots <a_{n-1}< a_n=\infty
\]
such that  for every $i=1,\cdots, n$  the interval $(a_{i-1}, a_i)$ contains precisely  one  critical value of $f$.

 Two  Morse functions $f,g: X\ra \bR$   will be called \emph{geometrically equivalent} if there exists an orientation preserving diffeomorphism $r:X\ra X$ and an orientation preserving diffeomorphism  $\ell: \bR\ra \bR$ such that
 \[
 g= \ell\circ f\circ r.
 \]
We denote by $\sim_g$ this   equivalence relation.

Two  Morse functions $f,g: X\ra \bR$   will be called \emph{topologically equivalent}    if  $\nu(f)=\nu(g)=n$ and  there exists  a  slicing $a_1<\cdots <a_{n-1}$ of $f$, a slicing  $b_1<\cdots <b_{n-1}$  of $g$  and orientation preserving diffeomorphisms
\[
\phi_i:\{f\leq a_i\}\ra \{g\leq b_i\},\;\;\forall i=1,\cdots , n,\;\;a_n=b_n=\infty.
\]
Two  Morse functions $f,g: X\ra \bR$  will be called \emph{homologically equivalent}  if  $\nu(f)=\nu(g)=n$ and  there exists  a  slicing $a_1<\cdots <a_{n-1}$ of $f$ and a slicing  $b_1<\cdots <b_{n-1}$  of $g$  such that for every $i=1,\cdots n$ the sublevel sets $\{f\leq a_i\}$ and $\{g\leq b_i\}$ have the same Betti numbers. We denote by $\sim_t$ and $\sim_h$ these equivalence  relations. Note that
\[
f_0\sim_g f_1 \Longrightarrow f_0\sim_t f_1\Longrightarrow  f_0\sim_h f_1.
\]
Set
\[
\m^n_X:=\{f\in \m_X;\;\;\nu(f)=n\}, \;\;\;[\m^n_X]_\ast:=\m^n_X/\sim_\ast,\;\;\ast\in\{h,t,g\}.
\]
Observe that we have natural projections
\[
[\m_X^n]_g\twoheadrightarrow [\m_X^n]_t\twoheadrightarrow [\m_X^n]_h.
\]
Since the excellent  Morse functions are stable  we deduce that  the geometric equivalence classes are open subsets of $\m_X^n$. This shows   that the quotient topology on $[\m_X^n]_\ast$, $\ast\in\{ g,t,h\}$,  is discrete.

If we think of a Morse function $f\in\m^n_X$ as defining a sort of ``triangulation'' on $X$ with $n$ simplices then  we can expect that $[\m_X^n]_g$ is finite\footnote{I haven't worked out a formal argument showing that $[\m_X]^n$ is finite, but I believe this to be the case.}  and therefore the sets $[\m_X^n]_{t,h}$ ought to be finite as well.   Clearly, the sets $\m_X^n$ contain information about the manifold $X$.  For example,   Morse theory implies that
\[
\m_X^n\neq \emptyset \Longrightarrow n\geq \sum_k b_k(X).
\]
We can then ask  how big are  of the sets $[\m^n_X]_\ast$, $\ast \in\{g,t,h\}$ and what kind of additional structure do they posses. It is more realistic to start by answering this question  for special $X$. The simplest case is $X=S^1$. We discuss briefly this case since it bares some similarities with the case $X=S^2$  we will discuss at length in this paper.

The sublevel sets of a Morse function  on $S^1$ are disjoint unions of closed intervals and we deduce that in this case the topological and homological classifications coincide.    A Morse function on $S^1$ has an even number of critical points,  but two of them are very special, namely the global minimum and maximum.

Suppose  the Morse function $f:S^1\ra \bR$ has $2m+2$ critical points $p_0, p_1, \cdots, p_{2m+1}$, with $c_i=f(p_i)<f(p_j)=c_j$ iff $i<j$.    Set $L_i:=\{f\leq c_i+\ve\}$. Then $L_0$ and $L_{2m}$ are closed intervals. If we set $x_i(f):= b_0(L_i)$ then we see that a Morse function  defines a  sequence
\[
\vec{x}(f):\{ 0,1,2,\cdots, 2m\}\ra \bZ_{>0},\;\; i\longmapsto{x}_i(f),
\]
satisfying
\[
x_0=x_{2m}=1,\;\; x_i>0,\;\;|x_{i+1}-x_i|=1,\;\;\forall i.
\]
Conversely, to any such sequence we can associate a  Morse function with $2m+2$ critical points  and we have
\[
f\sim_t g\Longleftrightarrow \vec{x}(f)=\vec{x}(g).
\]
We regard such a sequence  as a walk of length $2m$ on the lattice $\bZ\subset \bR$  with steps of size $\pm 1$  which starts and ends at $1$ and  it is confined to the positive chamber $\bZ_{>0}$. The number of  such walks can be easily determined using  Andr\'e's reflection principle, \cite[Ex. 14.8]{LW} and we obtain
\[
[\m^{2m+2}_{S^1}]_h= \binom{2m}{m}-\binom{2m}{m+2}=\frac{1}{m+2}\binom{2m+2}{m+1}=C_{m+1},
\]
where $C_k=\frac{1}{k+1}\binom{2k}{k}$ denotes the $k$-th Catalan number.

The number of geometric  equivalence classes of Morse functions on $S^1$ can be determined  using the calculus of snakes of Arnold, \cite{Ar,Ar1}.    We outline the main idea.

\begin{figure}[ht]
\centerline{\epsfig{figure=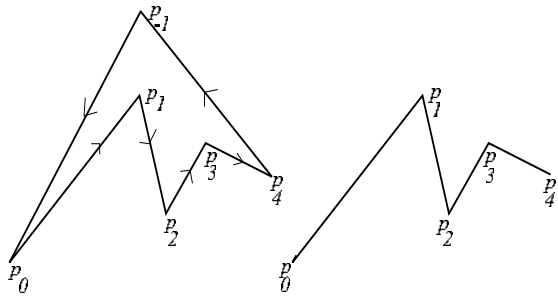,height=1.5in,width=2.8in}}
\caption{\sl Cutting a Morse function down to a snake.}
\label{fig: 0}
\end{figure}

Suppose $f: S^1\ra \bR$ is  Morse function with $2m+2$ critical points. Denote by $p_{-1}\in S^1$ the point where $f$ achieves its global maximum. Then, starting  at $p_{-1}$ label the critical points $p_0, \cdots, p_{2m}$, in counterclockwise  order (see Figure \ref{fig: 0}).  Now remove the arcs $[p_{-1},p_{0}), (p_{2m},p_{-1}]$ (see Figure \ref{fig: 0}).  What's left is  what Arnold calls a $A_{2m}$-snake. Their numbers are and the associated generating series are determined in \cite{Ar, Ar1}.

\begin{figure}[ht]
\centerline{\epsfig{figure=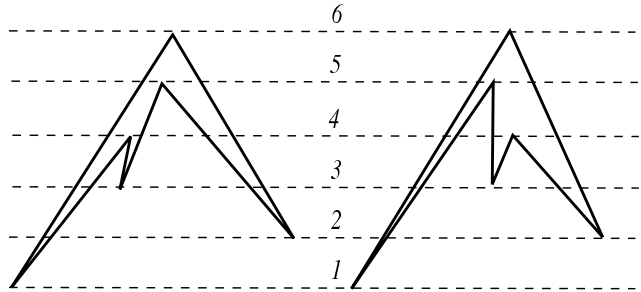,height=1.5in,width=3.3in}}
\caption{\sl Two Morse functions which are topologically equivalent but not geometrically equivalent.}
\label{fig: 01}
\end{figure}

Let us point out that the geometric classification of Morse functions on $S^1$ is different from the topological classification. For example the   Morse functions depicted in Figure \ref{fig: 01} are   topologically equivalent, yet not geometrically so. In this figure the Morse function is the height function $y$, the vertices are the critical points and the numbers attached to them are the corresponding critical values.

Note that   for the function  on the left hand side the first and third critical points lie in the same component of the sublevel set $\{y \leq 4+\ve\}$ but this is not the case for the function on the right hand side.

\section{The main results}

In this paper we  investigate   the structure  of the  sets $[\m^n_{S^2}]_\ast$ where $\ast\in\{h,g\}$.  Here is briefly its content.

In Section 3, using the basics of Morse theory we analyze what kind of handle additions can occur   as we cross the critical levels of  a  Morse function on  $S^2$.  This leads  immediately to a   bijection between the homology classes of such Morse functions and     sets of  lattice paths in $\bZ^2_{>0}$.    The cardinalities of such sets are  computed explicitly  in Theorem \ref{th: weak} using the reflection technique    of Gessel and Zeilberger, \cite{GZ}.  More precisely, in Theorem \ref{th: weak} we prove that the number of homology classes of Morse functions on $S^2$ with $2n$ critical points is
\[
h_{2n}=C_nC_{n-1},
\]
where $C_n$ denotes the $n$-th Catalan number.

In  Section 5 we describe a bijection between the set of geometric equivalence classes of Morse functions on $S^2$ and the set of \emph{Morse trees}. These are labelled trees  with vertices of degree one or  three such that   every node (i.e. vertex  of degree three) has at least one neighbor  with a bigger label and at least one neighbor with a smaller  label.  This correspondence   from  Morse functions to labelled graphs   first appeared on the mathematical scene in the work of G. Reeb \cite{Reeb1}.

A Morse tree has a natural   plane structure, i.e. a natural way of   linearly ordering the branches at every vertex.   In particular, every node will have a canonical pair of neighbors called \emph{successors}.  Moreover  the edges   can be canonically  colored with three colors  called  North, Northeast  and Southeast     such   every node has a unique   Northern  \emph{successor},  meaning that the edge connecting them  is a North edge. We called the resulting edge colored  plane tree a \emph{Morse profile} because the  tubular neighborhoods  of two Morse trees in $\bR^3$ ``look the same to the naked eye''  if they have the same profile.

In Section  6 we count the number of Morse profiles.  We show in Theorem \ref{th: profile} that their generating function  $f=f(t)$ satisfies the cubic equation $f(1-tf)^2=1$ and then we determine the  Taylor coefficients of this generating function using the Lagrange inversion formula.

In Section 7 we produce a two parameter recurrence for the number of Morse trees (Theorem \ref{th: main}) which  is computationally very  effective. In Section 8 we  associate to this  two parameter family of numbers  an  exponential type generating function of two independent variables and we prove that it satisfies  a singular  first order quasilinear p.d.e.   Fortunately,  the singularities of this equation can be resolved using a monoidal  blowup involving both the dependent and the independent variables. The new  equation is amenable to the classical method of characteristics \cite{CH} which yields a concrete description of the exponential  generating function of the number of geometric classes.

More precisely, if    $\xi_{n}$ denotes the number of   geometric classes of Morse functions with $2n+2$ critical points and
\[
\xi(t)= \sum_{n\geq 0}\xi_n\frac{t^{2n+1}}{(2n+1)!}
\]
then  in Theorem \ref{th: main1} we show that $\xi$ is the compositional inverse of the function
\[
\theta(s)=\int_0^s\frac{d\tau}{\sqrt{\frac{\tau^4}{4}-\tau^2+2s\tau+1}},
\]
i.e. $\xi(\theta(s))=s$. The proof also     reveals a relationship between   the Morse function count and the classical  Weierstrass elliptic functions $\wp$.

The Lagrange inversion formula can in principle be used  to compute the Taylor  coefficients of $\xi$, but  the  double recurrence in Theorem \ref{th: main} seems easier to implement  on a computer.

The set  of topological equivalence classes of Morse functions    seems very mysterious at this time. It has resisted all our attempts  to  uncover  a computationally  friendly  structure.   In section 9 present  some partial results.  Using   some results of Stanley \cite{St3}  on  counts  of paths  in the Young lattice of partitions  we were able to produce  a (far from optimal) lower bound.  There exists at least $1\cdots 3\cdots (2n+1)$ topological equivalence classes of   Morse functions on $S^2$ with  $2n+2$ critical points. We also compute the number of topological equivalence classes of Morse functions with at most $10$ critical points.

\section{The anatomy of a Morse function on the $2$-sphere}
\setcounter{equation}{0}

Denote by $[n]$ the  set $\{1,2,\cdots ,n\}$. Suppose $f$ is a excellent Morse function on the $2$-sphere. In this case $\nu(f)$ is a positive even integer and we set
 \[
 \bar{\nu}(f):=\frac{1}{2}\nu(f) -1.
 \]
We  say  that $\bar{\nu}(f)$ is the \emph{order} of $f$. We denote by  $\Cr_f\subset S^2$ the set of critical points and by $\Delta_f\subset \bR$ the set of critical values. Suppose $\bar{\nu}(f)=n$, i.e. $\#\Cr_f=2n+2$. The \emph{level function} associated to $f$ is the bijection
\[
\ell:\Cr_f\ra [2n+2],\;\; \ell (p):= \# \bigl\{\, q\in\Cr_f;\;\;f(q)\leq f(p)\,\bigr\}.
\]
For every $i\in[2m+2]$ we denote by $p_i$ the critical point of $f$ such that $\ell(p)=p_i$. Note that we have a bijection
\[
\hat{\ell}: \Delta_f\ra [2n+2],\;\; \Delta_f\stackrel{f^{-1}}{\Lra}\Cr_f\stackrel{\ell}{\Lra} [2n+2].
\]
We say that a critical value $c$ has level $i$ if $\hat{\ell}(c)=i$.

For every  regular value  $\min f < c <\max f$ of $f$    every
connected component of the sublevel  set $\{f\leq c\}$  is   a
sphere with  some disks removed. As we cross a critical value the
sublevel sets $\{f\leq c\}$      undergoes  a change  of type
$H_0$, $H_2$, $H_1^+$ or $H_1^-$.

\begin{figure}[ht]
\centerline{\epsfig{figure=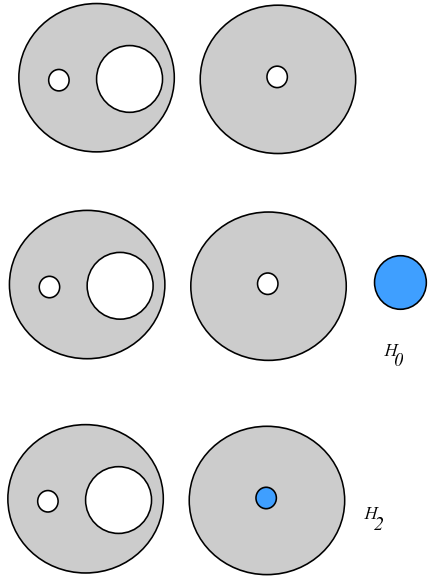,height=3in,width=2.2in}}
\caption{\sl Attaching a $0$-, $2$-handle.} \label{fig: 1}
\end{figure}

The change $H_0$ corresponds to crossing a local minimum of $f$.
The sublevel set    acquires a new component diffeomorphic to  a
closed disk (see Figure \ref{fig: 1} where the new component is
depicted in blue). The change $H_2$ corresponds to  crossing a
local maximum  of $f$  and   consists of attaching a $2$-disk to a
boundary component of $\{f\leq c\}$.

\begin{figure}[ht]
\centerline{\epsfig{figure=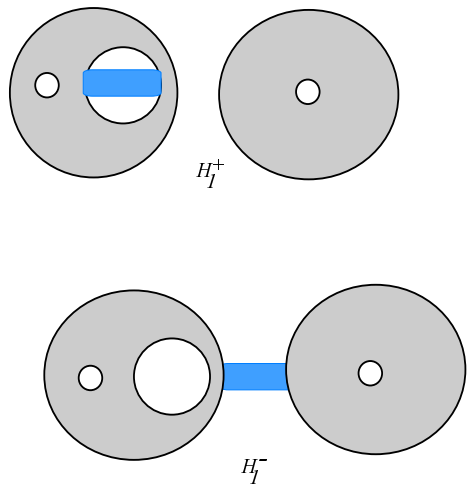,height=2.5in,width=2.4in}}
\caption{\sl Admissible  $1$-handle attachments.} \label{fig: 2}
\end{figure}

\begin{figure}[ht]
\centerline{\epsfig{figure=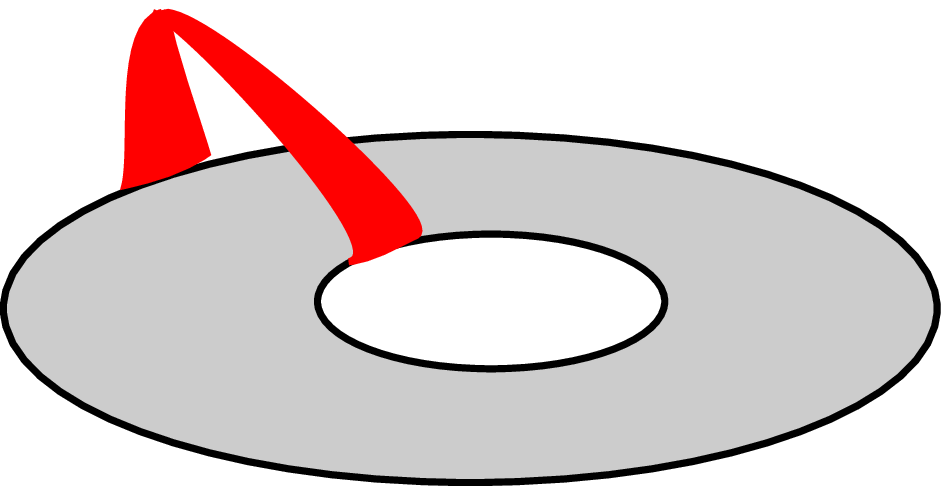,height=1.0in,width=1.7in}}
\caption{\sl Inadmissible   $1$-handle attachment.} \label{fig: 3}
\end{figure}

Crossing a saddle point  has the effect of attaching a $1$-handle.
This can be done in only two ways: $H_1^+$, $H_1^-$ (see  Figure
\ref{fig: 2}).  The case $H_1^+$ corresponds to the  gluing  of
the $1$-handle to the same boundary component of  the sublevel
set. The change $H_1^-$ corresponds to the attachment of the
$1$-handle to  different components of  the sublevel set.     A
priori, there is  another possibility, that of attaching the
$1$-handle to the same component  of $\{f\leq c\}$ but in
different boundary components. This  move can be excluded since it
produces   a $1$-cycle which cannot be  removed by     future
handle additions (see Figure \ref{fig: 3}).

To every regular value $\min f <c <\max f$ we associate  a  vector
\begin{equation}
\vec{b}=\vec{b}(c)= \bigl(\, b_0(X^c)\, ,\, b_1(X^c)+1\,\bigr),\;\;X^c=\{f\leq c\}.
\label{eq: betti}
\end{equation}
As we cross  a critical value the  vector $\vec{b}$
undergoes one of the changes  below
\[
\vec{b} \stackrel{H_0}{\Lra} \vec{b}+
(1,0),\;\;\vec{b}\stackrel{H_2}{\Lra} \vec{b}+ (0,-1),
\]
\[
\vec{b}\stackrel{H_1^+}{\Lra} \vec{b}+
(0,1),\;\;\vec{b}\stackrel{H_1^-}{\Lra} \vec{b}+ (-1,0).
\]
We denote by $m$ the number of \emph{local but not global} minima, by $M$ the number of \emph{local but not global} maxima, by $n_\pm$ the number of  $H_1^\pm$ handles. Note that these integers are  constrained by the    equalities
 \begin{equation}
 \left\{
 \begin{array}{rcl}
 m+n_-+n_++M &= &2\nu\\
 m-n_--n_++M& = & 0\\
 n_- &= & m\\
 n_+ &= & M
 \end{array}
 \right..
 \label{eq: min-nodes}
 \end{equation}
 The Morse polynomial of $f$ is therefore
 \[
 (1+m)+ \nu t+ (1+\nu-m)t^2.
 \]

\section{Counting homology equivalence classes}
\setcounter{equation}{0}

We denote by $h(n)$ the number of $\sim_h$-equivalence classes of Morse functions $f: S^2\ra \bR$ such that $\bar{\nu}(f)=n$.
\[
h(z)=\sum_{n\geq 0} h(n) z^n\in \bZ[[z]].
\]
We would like to compute  $h(z)$.

\begin{figure}[ht]
\centerline{\epsfig{figure=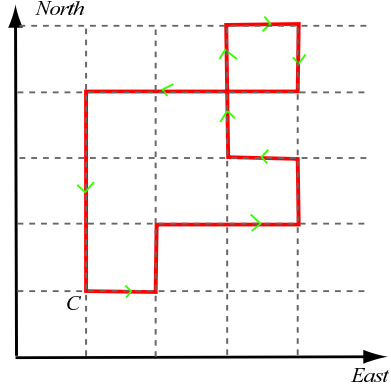,height=2.0in,width=2.0in}}
\caption{\sl A good path of type consisting of $4$ steps East and $4$ steps North.}
\label{fig: 4}
\end{figure}

Let $f\in \m_S^2$ and suppose its critical values are $c_0<c_1<\cdots <c_{2\bar{\nu}+1}$.  For $i=0,\cdots, 2\bar{\nu}$   and $0 <\ve  <\min(c_{i+1}-c_i)$ we set
\[
d_i:= c_i+\ve,\;\; 0\leq i\leq 2\bar{\nu}.
\]
The real numbers $(d_i)$ define a slicing of $f$.   To the
function  $f$ we now associate  using  (\ref{eq: betti}) the
sequence of lattice points
\[
P_i= \vec{b}(d_i)= (b_{0,i},b_{1,i}+1)\in \bZ_{>0}^2,\;\; b_{k,i}:= \dim H_k\bigl(\, \{f\leq d_i\},\,\bR\,\bigr).
\]
Set $e_1=(1,0)$, $e_2=(0,1)$. Observe that $P_0= P_{2\bar{\nu}}=(1,1)$ and  for every $i$ we have
\[
P_{i+1}-P_i\in \{\pm e_1, \pm e_2\}=:{\bf S}.
\]
We regard the sequence $\p(f)=\{P_0, P_1,\cdots, P_{2\bar{\nu}}\}$ as a  path  in the interior of the first quadrant such that  every step $P_i\ra P_{i+1}$  has length $1$  and is performed in one of the four possible   lattice directions at $P_i$,  East, North, West, South (see Figure \ref{fig: 4}).  Note that
\[
f\sim_h g \Longleftrightarrow \p_f=\p_g.
\]

Define now a \emph{lattice path} to be a  finite sequence of points $\gamma=\{P_0,\cdots, P_L\}$ such that $P_{i+1}-P_i\in {\bf S}$. The integer $L$ is called the \emph{length} of the path and it is denoted by $L(\gamma)$. The path is called \emph{good}  if  all the points $P_i$ are in the interior of the first quadrant. A path is called \emph{bad} if it is not good. For every integer $n$ we denote by $\p(P, Q; n)$ the set  of  paths of length $n$ starting at $P$ and ending at $Q$. $\p_g$ and $\p_b$ denote the subsets  consisting of good and respectively  bad paths. We set
\[
N(P,Q;n)=\#\p(P, Q;n),\;\;N_g(P, Q;n)=\#\p_g(P,Q;n),
\]
\[
N_b(P,Q;n)=\#\p_b(P, Q;n).
\]
The above discussion shows that
\[
h(n)=N_g(C,C;2n),\;\;C=(1,1).
\]
The number of good paths is computed in \cite{GZ}.  Consider two points  $P,Q$ in the first quadrant. Denote by $r_1, r_2:\bR^2\ra \bR^2$ the reflections
\[
r_1(t_1,t_2)=(-t_1,t_2),\;\;r_2(t_1,x_2)=(t_1,-t_2).
\]
These reflections generate the Klein group
\[
K=\bigl\{\,\one, r_1, r_2, r=r_1r_2\,\bigr\}\cong \bZ/2\times \bZ/2.
\]
Then as in the proof of \cite[Thm. 1]{GZ} we deduce
\[
N_b(P, Q;n)+ N_b(P, r(Q);n)= N_b(P,r_1(Q);n)+ N_b(P,r_2(Q); n).
\]
Since $r(Q),r_1(Q), r_2(Q)$  are not in the  first quadrant we deduce
\[
N_b(P,\rho(Q))=N(P,\rho(Q)),\;\; \forall \rho\in \{r,r_1,r_2\}
\]
and thus
\begin{equation}
N_g(P,Q;n)=N(P,Q;n)-N(P,r_1(Q);n)- N(P,r_2(Q); n)+ N(P, r(Q);n).
\label{eq: good-paths}
\end{equation}
Observing that
\[
N(P,Q; n)= N(0, Q-P; n)
\]
we deduce
\begin{equation}
N_g(C,C;n)=N(0,0;n)-N(0,-2e_2;n)- N(0,-2e_1; n)+ N(0, -2(e_1+e_2);n).
\label{eq: good-paths2}
\end{equation}
For $t=(t_1,t_2)$ and $\vec{x}=(x_1,x_2)\in\bZ^2$ we set $t^{\vec{x}}=t_1^{x_1}t_2^{x_2}$. Consider the \emph{step} polynomial
\[
\s(t):=\sum_{\vec{x}\in {\bf S}}t^{\vec{x}}=t_1+t_2+t_1^{-1}+t_2^{-1} \in \bZ[[t_1,t_2,t_1^{-1},t_2^{-1}]].
\]
For every $\vec{x}\in \bZ^2$ define
\[
\C_{\vec{x}}: \bC[[t_1,t_2,t_1^{-1},t_2^{-1}]]\ra \bC, \;\;\;A=\sum_{\vec{y}\in\bZ^2} A_{\vec{y}}t^{\vec{y}}\longmapsto A_{\vec{x}}.
\]
Equivalently, if we regard $t_j$ as a complex parameter, $t_j=|t_j| e^{\ii\theta_j}$, we  have the integral formula
\[
\C_{\vec{x}}(A)=\oint A(t)t^{-\vec{x}}  :=\frac{1}{4\pi^2}\int_{|t_1|=|t_2|=1} A(t)t^{-\vec{x}} d\theta_1d\theta_2.
\]
For every $\vec{x}\in \bZ^2$ we then have
\[
N(0,\vec{x};n)=\C_{\vec{x}}\bigl(\, \s(t)^n\,\bigr),
\]
so that
\[
N_g(C,C;n)=\oint \s(t)^n(1-t_1^2-t_2^2+t_1^2t_2^2) = \oint\s(t)^n(t_1^2-1)(t_2^2-1) .
\]
Hence
\begin{equation}
h(z^2)= \oint\frac{(t_1^2-1)(t_2^2-1)}{1-z\s} dt=\oint\underbrace{\frac{(t_1^2-1)(t_2^2-1)}{1-z\bigl(\, t_1+t_1^{-1}+t_2+t_2^{-1}\,\bigr)}}_{=:K(z,t)} dt.
\label{eq: weak-morse}
\end{equation}
Observe that
\[
K(z,t)= \frac{t_1t_2(t_1^2-1)(t_2^2-1)}{t_1t_2-z(t_1t_2+1)(t_1+t_2)}\;\;{\rm and}\;\;K(z,t_1,t_2)=K(z,t_2,t_1).
\]
 We write
\[
u=t_1+t_2,\;\; v=t_1t_2
\]
and we deduce
\[
K(z,t)= \frac{v(v^2-u^2+2v+1)}{v-zu(v+1)}= (v^2-u^2+2v+1)\cdot \frac{1}{1-zu(1+1/v)}
\]
\[
=(v^2-u^2+2v+1)\cdot\sum_{n\geq 0}z^nu^n\bigl(\,1+v^{-1}\,\bigr)^n
\]
\[
=\sum_{n\geq 0} z^n\sum_{k=0}^n \binom{n}{k} (v^2-u^2+2v+1)u^nv^{-k}
\]
\[
=\sum_{n\geq 0} z^n\sum_{k=0}^n\binom{n}{k}(u^nv^{2-k}-u^{n+2}v^{-k}+2u^nv^{1-k}+u^nv^{-k}).
\]
Fortunately, very few terms in this sum contribute to (\ref{eq: weak-morse}) since
\[
\oint u^nv^{-k}=0,\;\;\forall n\neq 2k.
\]
Moreover
\[
\oint u^{2m}v^{-m}= \binom{2m}{m}.
\]
We deduce
\[
\oint K(z,t)
\]
\[
=\sum_{m\geq 0} z^{2m}\Biggl(\sum_{k=0}^{2m}\binom{2m}{k}\Bigl(\, \oint u^{2m}v^{2-k}-\oint u^{2m+2}v^{-k}+2\oint u^{2m}v^{1-k}+\oint u^{2m}v^{-k}\,\Bigr)\,\Biggr).
\]
This shows that for every $m\geq 0$ we have
\[
h(m)=  \underbrace{\binom{2m}{m+2}\binom{2m}{m}}_{a}-\underbrace{\binom{2m+2}{m+1}\binom{2m}{m+1}}_{b}+2\;\underbrace{\binom{2m}{m}\binom{2m}{m+1}}_{c}+\underbrace{\binom{2m}{m}\binom{2m}{m}}_{d}.
\]
Now observe that
\[
a+c=\binom{2m}{m+2}\binom{2m}{m}+ \binom{2m}{m}\binom{2m}{m+1}=\binom{2m}{m}\binom{2m+1}{m+2},
\]
\[
c+d=\binom{2m}{m}\binom{2m}{m+1}+ \binom{2m}{m}\binom{2m}{m}=\binom{2m}{m}\binom{2m+1}{m+1}.
\]
The sum of the right-hand-sides of the above equalities is
\[
a+2c+d=\binom{2m}{m}\binom{2m+1}{m+2}+ \binom{2m}{m}\binom{2m+1}{m+1}=\binom{2m}{m}\binom{2m+2}{m+2}.
\]
Hence
\[
h(m)=a+2c+d-b=\binom{2m}{m}\binom{2m+2}{m+2}-\binom{2m+2}{m+1}\binom{2m}{m+1}
\]
\[
= \frac{1}{(m+2)(m+1}\binom{2m+2}{m+1}\binom{2m}{m}=C_{m+1}C_{m},
\]
where $C_m=\frac{1}{m+1}\binom{2m}{m}$ denotes the $m$-th  Catalan number. We have thus proved the following result.

\begin{theorem}
\[
h(m)=C_{m+1}C_m,\;\;\forall m\geq 0.
\]
\label{th: weak}
\end{theorem}

\begin{ex}   Let us test the validity of the above formula   for small values of $m$.   We have
\[
C_0=1,\;\;C_1=1,\; C_2=2,\;\;C_3=5,\;\;C_4=14,\;\;{\rm etc.}
\]
We deduce
\[
h(1)=2,\;\;h(2)=10.
\]
The two  good paths of length $2$ from $C$ to $C$ are
\[
\{E,W\},\;\;\{N,S\}.
\]
To list the  ten good  paths of length $4$  it suffices to list
only the five  good paths   whose first step  points East. The
other five are obtained  via  the reflection in the diagonal
$x_1=x_2$ given by
\[
E\longleftrightarrow N,\;\; W\longleftrightarrow S.
\]
The five good path with initial Eastbound step are
\[
\{E,W,E,W\},\;\;\{ E,E,W,W\},\;\;\{E,N,S,W\},\;\; \{E,W, N,S\},\;\;\{E,N,W,S\}.\proofend
\]
\end{ex}

\section{Invariants of geometric equivalence classes}
\setcounter{equation}{0}

 In this section we describe the main  combinatorial invariants  associated to a  geometric equivalence class of Morse functions.

 Suppose $f$ is a   Morse function $f: S^2\ra \bR$ of order $n$ with critical points $\{p_1,\cdots, p_{2n+2}\}$ labelled such that
 \[
 f(p_i)<f(p_j)\Longleftrightarrow i< j.
 \]
In the sequel, by a \emph{labelling} of a graph will understand an injection from the set of vertices to the set of real numbers. For every graph $\Gamma$ and every  labelling  $\vfi$ of $\Gamma$  we define the \emph{level function} associated to $\vfi$
\[
\ell_\vfi:\cv\ra   \bZ,\;\; \ell_\vfi(v):=\#\bigl\{\, u\in \cv(\Gamma);\;\;\vfi(u)\leq \vfi(v)\,\}.\proofend
\]
We will  construct inductively a sequence of labelled forests  $F_1, \cdots,  F_{2n+1}, F_{2n+2}$ which replicates the changes in the  level sets $L_k:=\{f \leq c_k+\ve\,\}$, $\ve <\min_k(c_{k+1}-c_k)$.

The vertices of $F_k$ are   colored in black  and white. The black vertices of $F_k$ are bijectively labelled  with  the labels $\{1, 2,\cdots,k\}$.  The  white   vertices of $F_k$ have degree one and correspond to the  boundary components of the level set $L_k$ which we know is a disjoint union of  spheres with disks removed.

\begin{figure}[ht]
\centerline{\epsfig{figure=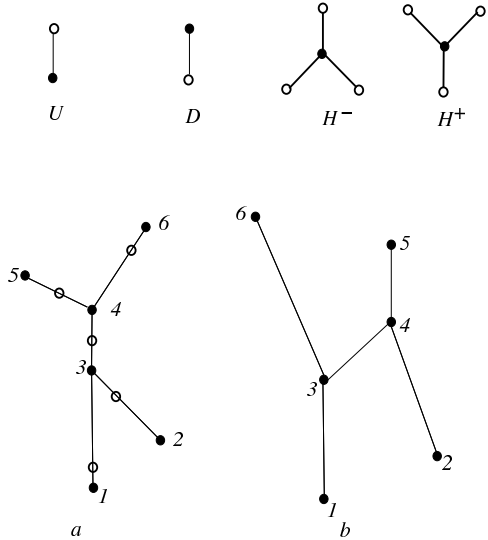,height=2.8in,width=2.5in}}
\caption{\sl Simulating handle attachments.}
\label{fig: 5}
\end{figure}

The construction of this  sequence of forests goes as follows.  $F_1$ is the  \emph{up} tree  $U$ depicted in Figure \ref{fig: 5} with the black vertex   labelled $1$.

The passage    $F_k\ra F_{k+1}$  depends  on the topological   change $L_k\ra L_{k+1}$, i.e. on the nature of the handle attachment   as we cross the critical value of level $k+1$. As explained in  Section 3, there are four types of handle attachments $H_0,H_2, H_1^\pm$.

\smallskip

\noindent $\bullet$ $H_0$-change.     $F_{k+1}= F_k \sqcup U$. The new black vertex is labelled  $k+1$.  The labels of the previous black vertices are not changed.

\noindent $\bullet$  $H_2$-change.  Glue the white vertex of the \emph{down} tree $D$ in Figure \ref{fig: 5} to a  white vertex of $F_k$ corresponding to the component of $\pa L_k$ where this two handle is attached. The new black vertex is labelled  $k+1$.

\noindent $\bullet$  $H^-_1$-change.    Glue the two lower white vertices of the tree $H^-$ in Figure \ref{fig: 5} to two   white vertices of $F_k$ situated in different path components. Label the new  black vertex $k+1$ and ignore the two white vertices where the gluing took place.

\noindent $\bullet$   $H_1^+$-change. Glue the lower white vertex of $H^-$ to the white vertex of $F_k$ corresponding to the  component of $\pa L_k$ where this $1$-handle is attached. Label the new  black vertex $k+1$ and ignore the  white vertex where the gluing took place.

\smallskip

The last forest, $F_{2n+2}$ consists of single a tree.   We will refer to it as the \emph{labelled tree} associated to $f$ and we will denote it by $\Gamma_f$. We denote by $\ell_f:\Gamma_f\ra [2n+2]$ the associated labelling.

We say that a vertex $v$ is higher than a vertex $u$ (or that $u$ is lower than $v$) if $\ell_f(v)>\ell_f(u)$. This  tree  has $2n+2$ vertices, and each of them has either one or three neighbors.  We will refer to the vertices of degree $3$ as \emph{nodes}.   Note that each node  has a higher  neighbor and a lower  neighbor.

To reconstruct the functions from the labelled tree   represent it as a tree  embedded   in $\bR^3$ so that the following hold.

\smallskip

\noindent $\bullet$ The $z$ coordinates of the  vertices are equal to the labels.

\noindent $\bullet$ The edges are smoothly embedded arcs.

\noindent $\bullet$  The restriction of the height function  $z$ to any edge is an increasing function, or equivalently, the  arcs representing the edges have no horizontal tangents.

\smallskip

We identify $S^2$ with the boundary of a thin tubular neighborhood of the tree in $\bR^3$ (imagine the edges made   of thicker and thicker spaghetti). The Morse function is then the restriction of the  height function $z$ to  this surface.

\begin{definition} (a) A \emph{normalized Morse tree  order $n$} is a pair $(\Gamma,\ell)$, where $\Gamma$  is a tree  with vertex set $\cv(\Gamma)$ of cardinality  $2n+2$  and $\ell: \cv(\Gamma)\ra [2n+2]$ is a labelling such that the following hold.

\smallskip

\noindent (a1) Every vertex of $\Gamma$ has degree one or three.  We will refer to  the degree $3$ vertices as \emph{nodes}.

\noindent (a2) Every node has   at least one lower neighbor and at least one higher neighbor.

\smallskip

 A degree $1$ vertex is called a \emph{maximum/minimum} if it is higher/lower than its unique neighbor. We denote by $\n(\Gamma)\subset\cv(\Gamma)$ the set nodes,  and by $\m^\pm(\Gamma)$ the set of  maxima/minima.

\smallskip

\noindent (b) Two normalized Morse  trees $(\Gamma_i,\ell_i)$, $i=0,1$ are said to be \emph{isomorphic}, $(\Gamma_0,\ell_0)\cong (\Gamma_1,\ell_1)$,  if there exists a bijection
\[
\phi: \cv(\Gamma_0)\ra \cv(\Gamma_1)
\]
such that $u,v\in \cv(\Gamma_0)$ are neighbors in $\Gamma_0$ iff $\phi(u),\phi(v)$ are  neighbors in $\Gamma_1$ and
\[
\ell_1(\phi(u))=\ell_0(u),\;\;\forall u\in\Gamma_0. \proofend
\]
\end{definition}

Observe that the order of a (normalized) Morse tree  is precisely the number of nodes.   The  labelled tree  $(\Gamma_f,\ell_f)$  we have associated to  a Morse function $f$ on $S^2$ is a  normalized Morse tree and we see that
\begin{equation}
f\sim_g g \Longleftrightarrow (\Gamma_f,\ell_f)\cong (\Gamma_g,\ell_g).
\label{eq: morse-tree}
\end{equation}
We denote by $\e(\Gamma)$ the set of edges.  Since $\Gamma$ is a tree  with $2n+2$ vertices we deduce
\[
\#\e(\Gamma)= 2n+1.
\]

\begin{definition} A \emph{Morse tree} is a pair $(\Gamma,\vfi)$, where $\Gamma$ is a tree  and  $\vfi$ is a labelling of $\Gamma$ such that $(\Gamma,\ell_\vfi)$ is a normalized Morse tree.
\qed
\end{definition}

Note that the set of Morse trees  is equipped with a natural involution
\[
(\Gamma,\vfi)\longleftrightarrow (\Gamma, -\vfi).
\]
We will refer to this as the \emph{Poincar\'e duality}.  We will also denote the Poincar\'e dual of a Morse tree $\Gamma$ by $\check{\Gamma}$. A Morse tree  will be called \emph{self-dual} if it is isomorphic to its Poincar\'e dual.

\begin{remark} (a)     The  above  labelled tree is sometime known as the \emph{Reeb graph} associated to the Morse function  $f$, \cite{BF, Ku, Reeb1, Sha}.  In \cite{Reeb1} Reeb associates a graph  to an excellent  Morse function $f: M\ra \bR$, $M$ smooth compact manifold. More precisely one can define an equivalence relation on $M$  by declaring two points $x_0,x_1\in M$ equivalent if and only if  they lie  in the same path component of some level set $f^{-1}(c)$. The   quotient space  equipped with the quotient topology is then homeomorphic with a finite one dimensional simplicial complex whose vertices correspond to the critical points of $f$. The vertices are then naturally labelled by the value of $f$ at the corresponding critical point. When $f$ is an excellent Morse function   on $S^2$ then its Reeb graph coincides with   the Morse tree  defined above.

 (b) It is very easy to produce  weaker invariants of a geometric equivalence class.  Given $u,v\in \m^-(\Gamma)$ there exists a unique path in $\Gamma_f$ connecting the vertices $u,v$. $m_f(u,v)$ is the level of the highest vertex along this path. $m_f(u,v)$ is the Mountain Pass level of $u$ and $v$.  In Figure \ref{fig: 5} we have depicted two  Morse functions $a,b$ which are  topologically equivalent. They are not geometrically  equivalent because
\[
m_a(1,2)= 3,\;\;m_b(1,2)=4.
\]
Given a  Morse tree $(\Gamma,\vfi)$   we a have function
\[
\delta_\vfi:V(\Gamma)\ra \cv(\Gamma),\;\;\delta_\vfi(v)=\mbox{lowest neighbor of $v$}.
\]

\begin{figure}[ht]
\centerline{\epsfig{figure=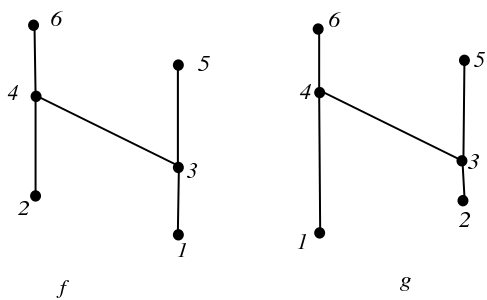,height=1.5in,width=2.5in}}
\caption{\sl Two Morse functions which are not geometrically equivalent, yet topologically equivalent.}
\label{fig: 6}
\end{figure}

In Figure \ref{fig: 6} we have depicted two  topologically equivalent Morse functions which are not geometrically  equivalent because $\delta_f(1)=3\neq 4=\delta_g(1)$. Alternatively, the critical points $\{1,3\}$ of $f$  lie in the same path component of $\{f\leq 3+\ve\}$  but the corresponding points of $g$  do not lie  in the same component of $\{g\leq 3+\ve\}$. Note however that $m_f(1,2)=m_g(1,2)=4$.
\qed
\end{remark}

\section{Planar Morse profiles}
\setcounter{equation}{0}

The only graphical representations of    Morse trees we have depicted so far in Figure \ref{fig: 5}, \ref{fig: 6}  are very special.  They are \emph{planar, cartesian}, i.e. the graphs are embedded in the cartesian plane $\bR^2$ (with coordinates $(x,y)$),  the labelling is given by the   height function $y$ and the edges  are  smooths  arcs in the plane such that the  tangent lines along them are  never horizontal\footnote{This means that the restriction of the height function along any arc has no critical point.}.

\begin{figure}[ht]
\centerline{\epsfig{figure=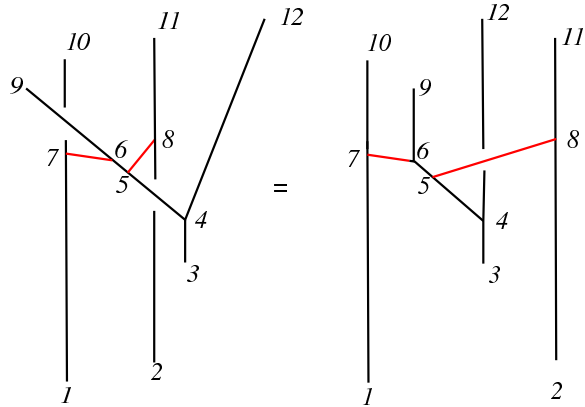,height=2.1in,width=3in}}
\caption{\sl A non planar Morse function?}
\label{fig: 14}
\end{figure}

Not  every Morse tree admits such a nice representation. For example, the Morse tree depicted in Figure \ref{fig: 14} does not seem to admit such a representation.

To understand the origin of the tree in  Figure \ref{fig: 14}  let us observe  that if we remove  the edges $[6,7]$ and $[5,8]$ then we obtain a collection of Morse trees  with no $H^-_1$-handles. These trees  are  determined by the    ``gradient'' flow.   We define the gradient  of a  Morse tree $(\Gamma,\vfi)$ to be the map
\[
\nabla: \cv(\Gamma)\ra \cv(\Gamma)
\]
where $\nabla(v)=v$ if $v$ is a local minimum while  if $v$ is a node then $\nabla(v)$  is the lowest neighbor of $v$  which lies below $v$.  Thus  if $v$ is not a local minimum  then the edge $[v,\nabla(v)]$ is the fastest descent edge at $v$.  The  gradient flow is    the (semi-)flow on $\cv(\Gamma)$ defined by the iterates of $\nabla$. Every orbit  ends at a local minimum and   the set of orbits ending at the same local minimum  determines a subtree  like in Figure \ref{fig: 14}.

However, every Morse tree admits  a canonical structure of plane tree.

\begin{proposition} Every normalized Morse tree $(\Gamma,\ell)$ of order  $n$ has a natural structure of  planted trivalent plane tree. i.e.  a rooted plane tree such  that the root has  degree $1$, and all other vertices have degrees $1$ or $3$. with $2n+2$ vertices. We denote it by $[\Gamma,\ell]$.
\end{proposition}

\noindent {\bf Proof}\hspace{.3cm}  Indeed, the root corresponds to the global minimum. There is only one edge starting at the root.   Label this edge by $1$.    Inductively,   once  we reached a vertex, then the unlabelled edges at this vertex can be ordered by the value  of the Morse function at their other endpoints.
\qed

\begin{figure}[ht]
\centerline{\epsfig{figure=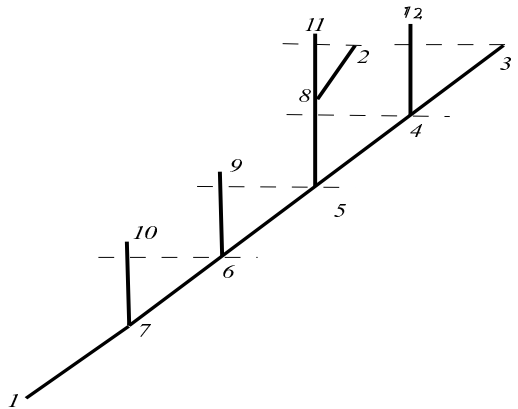,height=2.1in,width=2.7in}}
\caption{\sl A plane binary tree associated to a Morse function.}
\label{fig: 16}
\end{figure}

 As explained in \cite[Exer. 6.19.f]{St2} the number of such trees with $2n+2$ vertices is the Catalan number $C_n=\frac{1}{n+1}\binom{2n}{n}$. For example the planted ternary plane tree associated to the Morse tree in Figure \ref{fig: 14} is depicted in Figure \ref{fig: 16}.

We want to associate a more refined structure to the Morse tree, called the \emph{planar Morse profile}. This is obtained by drawing the associated plane tree in a very special way.

Every node $v$ of this plane tree  has two successors, a left successor, $l(v)$ and a right  successor, $r(v)$, where the left successor   of $v$ is defined to be the successor with highest label.  We will draw   the   edge connecting  $v$ to $l(v)$   by a vertical arrow pointing North. The  edge  $[v,r(v)]$ will be represented  by an arrow  pointing  Northeast  or Southeast. The arrow will point  Northeast  if the level  of $r(v)$ if  greater that the level of $v$, and will point  Southeast if the level of $r(v)$ is smaller than $v$.  The only  edge at the root will be vertical, and will point upward. After we do this then we remove the labels indicating the levels.

Formally, we are  coloring the edges withy three colors, North, Northeast and Southeast so that each node $u$ has exactly one  northern successor $v$, i.e. the edge $[u,v]$ is colored North.     We introduce  a partial order on the set of vertices  by declaring $u<v$ if the unique path from $u$ to $v$ consists of North or Northeast edges. This order has the property that every node has at least one   greater neighbor and at least one smaller neighbor.

\begin{ex} The planar profile of the tree depicted in Figure \ref{fig: 16}  is described in Figure \ref{fig: 17}.  For the reader's convenience we have not removed the labels so the correspondence with Figure \ref{fig: 16} is  more visible.

\begin{figure}[ht]
\centerline{\epsfig{figure=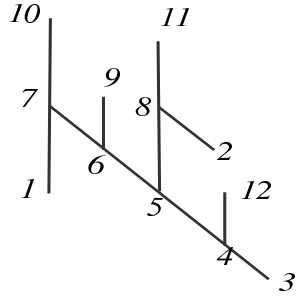,height=1.5in,width=1.5in}}
\caption{\sl A planar  Morse profile.}
\label{fig: 17}
\end{figure}

The Morse trees  depicted in Figure \ref{fig: 5} have different profiles, yet   the corresponding Morse functions are topologically equivalent. If Figure \ref{fig: 18} we have depicted the planar Morse profiles with $2$, $4$ and $6$ vertices.\qed

\begin{figure}[ht]
\centerline{\epsfig{figure=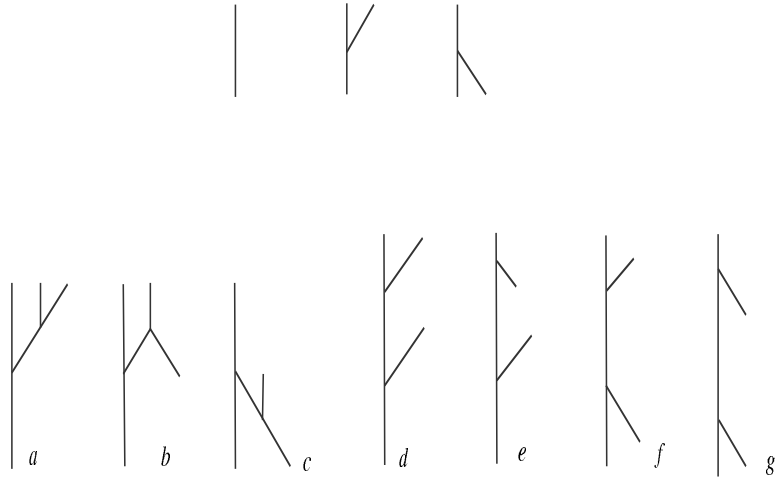,height=2.5in,width=4in}}
\caption{\sl The planar  Morse profiles with at most six vertices.}
\label{fig: 18}
\end{figure}

\end{ex}

Denote by $f_n$ the number of planar Morse profiles with $2n+2$ vertices.  Figure \ref{fig: 18} shows that
\[
f_0=1, \;\;f_1=2,\;\;f_2=7.
\]
Set
\[
f(t)=\sum_{n\geq 0} f_nt^n.
\]
We would like to determine $f(t)$. To achieve this we need to investigate the anatomy of a Morse profile $\Gamma$.

First we denote by $\nu(\Gamma)$ the number of vertices and we set
\[
\bar{\nu}(\Gamma):=\frac{\nu-2}{2}\Longleftrightarrow \nu=2\bar{\nu}+2.
\]
 We denote by $\rho=\rho_\Gamma$ its root and by $\rho^+=\rho_\Gamma^+$ its only neighbor.  The \emph{right wing} of the Morse profile is the maximal sequence of vertices $\{v_0,\cdots, v_s\}$, where
\[
v_0=\rho_\Gamma^+, \;\; v_{i+1}=r(v_i),\;\;i=0,\cdots, s-1.
\]
The maximality assumption implies that the last vertex $v_s$ on the right wing is a degree $1$ vertex of the tree.

A vertex  $v_i$ on the   right wing is called a \emph{turning point} if  either $v_i$  is the last point on the path, or $v_0$ is the first point on the path and $v_0>v_1$ or $0<i<s$ and $v_{i-1}< v_i >v_{i+1}$.

\begin{lemma}[Monotonicity] If $\{v_0,\cdots,v_s\}$ is  the right wing of $\Gamma$ and $v_{i-1},v_i,v_{i+1}$  are  three consecutive vertices along this path such that $v_{i-1}> v_i$ then $v_i>v_{i+1}$.
\end{lemma}
\noindent {\bf Proof}\hspace{.3cm} Indeed the edge $[v_i, v_{i+1}]$ cannot point Northeast    since in that case $v_i$ would have three  higher  neighbors $v_{i-1},l(v_i), v_{i+1}$. \qed

\begin{figure}[ht]
\centerline{\epsfig{figure=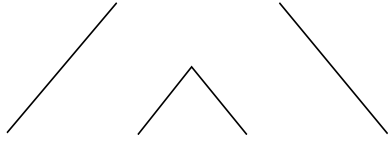,height=0.7in,width=2in}}
\caption{\sl The possible shapes of the right wing.}
\label{fig: 19}
\end{figure}

The monotonicity lemma implies that there exists exactly one turning point on the right wing, so that it  must have one of the three shapes depicted in Figure \ref{fig: 19}. In other words, the right wing can  change direction at most once.

\begin{figure}[ht]
\centerline{\epsfig{figure=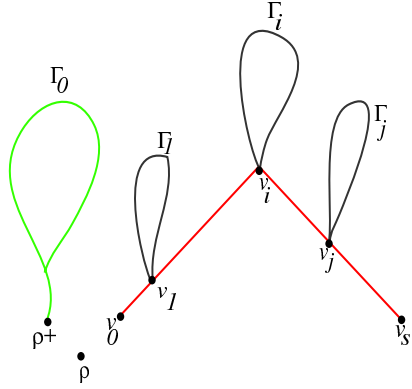,height=2in,width=2.1in}}
\caption{\sl Deconstructing a Morse profile.}
\label{fig: 20}
\end{figure}

Suppose that we now cut a Morse profile at the vertex $\rho_\Gamma^+$.  We  obtain  a configuration of the type depicted in Figure \ref{fig: 20}, where $\{v_0,\cdots, v_s\}$ denotes the right wing, $\Gamma_j$ denotes the left branch of $\Gamma$ at $v_j$, and $v_i$ is the  turning point of the right wing.

The Morse profile is completely determined by the following data.

\smallskip

\noindent $\bullet$ The the length $s$ of the right wing.

\noindent $\bullet$ The location $i\in\{0, 1,\cdots, s\}$ of  the turning point where the  right wing turns\footnote{If it never turns Southeast then we declare the turning point to be the last vertex along the right wing, that is $i=s$.} Southeast.

\noindent $\bullet$ The Morse profiles $\Gamma_0,\cdots, \Gamma_{s-1}$ satisfying
\[
\nu(\Gamma_0)+\cdots \nu(\Gamma_{s-1})= \nu(\Gamma)-2\Longleftrightarrow\bar{\nu}(\Gamma_0)+\cdots +\bar{\nu}(\Gamma_{s-1})=\bar{\nu}(\Gamma)-s.
\]
Observe that since $\deg(v_s)=1$ the tree $\Gamma_s$ is empty.

\smallskip

Since the  turning point can be chosen in $(s+1)$-different ways   we deduce
\[
f(n) =\sum_{s=1}^n\Bigl( \sum_{n_0+\cdots + n_{s-1}=n-s} (s+1) f_{n_0}\cdots f_{n_{s-1}}\Bigr).
\]
Equivalently, for every $n\geq 1$ we have
\[
f_nt^n=\sum_{s=1}^n(s+1)t^s\Bigl( \sum_{n_0+\cdots + n_{s-1}=n-s} f_{n_0}\cdots f_{n_{s-1}})t^{n_0+\cdots n_{s-1}}\,\Bigr).
\]
If we sum over $n$ we deduce
\[
f(t)-f_0= \sum_{s\geq 1}  (s+1) (tf)^s\Longleftrightarrow f=\sum_{s\geq 0} (s+1) (tf)^s=\frac{1}{(1-tf)^2}.
\]
It is more convenient to rewrite the last equality as
\[
tf=\frac{t}{(1-tf)^2}.
\]
Thus, if we set $w=w(t)=tf(t)$ and $Z= w(1-w)^2$ we deduce $Z(w(t))=t$  so that $w$ is the compositional inverse of $Z$.  Using the Lagrange inversion formula \cite[Thm. 5.4.2]{St2}  we deduce
\[
w=\sum_{n\geq 0} \frac{w_n}{n!}t^n,\;\;\; w_n= \frac{d^{n-1}}{dw^{n-1}}|_{w=0}\Bigl(\frac{w}{Z(w)}\Bigr)^n .
\]
Now observe that
\[
\frac{w}{Z(w)} =\frac{1}{(1-w)^2}\Longrightarrow \Bigl(\frac{w}{Z(w)}\Bigr)^n=\frac{1}{(1-w)^{2n}}.
\]
The coefficient of $w^{n-1}$ in $(1-w)^{-2n}$ is
\[
\frac{(2n)(2n+1)\cdots (3n-2)}{(n-1)!}
\]
and we deduce
\[
\frac{d^{n-1}}{dw^{n-1}}|_{w=0}(1-w)^{-2n}=(2n)(2n+1)\cdots (3n-2)
\]
so that
\[
f_{n-1}=\frac{w_n}{n!}=\frac{(2n)(2n+1)\cdots (3n-2)}{n!}=\frac{1}{n}\binom{3n-2}{n-1}\Longrightarrow f_n=\frac{1}{n+1} \binom{3n+1}{n}.
\]
Thus
\[
f_0= 1,\;\;f_1= \frac{1}{2}\binom{4}{1}= 2,\;\;f_2=\frac{1}{3}\binom{7}{2}=7,
\]
as   we discovered before. We have thus proved  the following result.

\begin{theorem} The number of planar Morse profiles with $2n+2$ vertices is
\[
f_n=\frac{1}{n+1} \binom{3n+1}{n}.
\]
Moreover, the generating series $f(t)=\sum_{n\geq 0} f_nt^n$ is a solution of the cubic equation
\[
f(1-tf)^2=1.\proofend
\]
\label{th: profile}
\end{theorem}

\section{A recursive construction of geometric equivalence classes}
\setcounter{equation}{0}

We would like to produce a two-parameter recursion formula for the number of normalized Morse trees of order $n$. One of the parameters will be the number of nodes and the  other parameter will be the level of the lowest node.

 If $(\Gamma,\vfi)$ is a  Morse tree and $\ell$ is the associated level function we set
\[
\lambda(\Gamma):= \min \{ \ell(v);\;\;v\in\n(\Gamma)\},\;\;\mu(\Gamma):=\lambda(\Gamma)-1.
\]
Observe that $\lambda(\Gamma)>1$ since $1$ is a global minimum of $\ell$ so that  $1\in \m^-(f)$. Hence $\mu(\Gamma)\geq 1$. If the order of $\Gamma$ is $n$ then the cardinality of $\n(\Gamma)$ is $n$. Moreover $2n+2\in E^+(f)$ and we deduce
\[
\lambda(\Gamma)\leq  n+2\Longrightarrow \mu(\Gamma)\leq n+1.
\]
In other words,   $\lambda(\Gamma)$ is the lowest level of a node. In particular we deduce that the vertices $1,2,\cdots,\mu(\Gamma)\}$ must be local minima.
This shows that if  $(\Gamma,\vfi)$ is an un-normalized Morse tree then $\mu(\Gamma)=m$ if an only if the first $m$ vertices of $\Gamma$ are local minima of $\vfi$ and the  while the vertex of  level $(m+1)$ is a saddle point.

Denote by $F_n(m)$ the number of  normalized Morse trees $(\Gamma,\ell)$ of order $n$ such that $\mu(\Gamma)\geq m$, that is Morse trees such that the lowest $m$ vertices are local minima.
Set
\[
f_n(m):= F_n(m)-F_n(m+1).
\]
Note that $f_n(m)$ denotes the number of Morse trees such that  the lowest  $m$ vertices are local minima while the $(m+1)$-th vertex is a node. Observe that $F_n(1)$ is precisely the number of  normalized Morse trees of order $n$ and
\[
F_n(m)=\sum_{k\geq m} f_n(k),\;\;  f_n(m)=0,\;\;\forall m>n+1.
\]
Since  a Morse tree of order $n$  has exactly $n+2$ extrema and one of them must be the highest vertex we deduce
\[
 F_n(m)=0,\;\;\forall  m>n+1\Longrightarrow f_n(n+1)=F_n(n+1).
 \]
\begin{theorem}
\[
f_n(m)= \binom{m}{2}F_{n-1}(m-1)
\]
\[
+\sum_{m_1+m_2=m-1}\Biggl(\,\sum_{n_1+n_2=n-1} \frac{m}{2}\binom{m-1}{m_1}\binom{2n-m+1}{2n_1-m_1+1}F_{n_1}(m_1+1) F_{n_2}(m_2+1)\,\Biggr).
\]
In particular, when $m=n+1$ we have
\[
F_n(n+1)=f_n(n+1)= \binom{n+1}{2} F_{n-1}(n).
\]
\label{th: main}
\end{theorem}

\noindent {\bf Proof}\hspace{.3cm}   For every  subset $C\subset \bR$ we   denote by $\f_C$ the set of un-normalized Morse trees $(\Gamma,\vfi)$ such that
\[
\vfi\bigl(\, \cv(\Gamma)\,\bigr)=C.
\]
We say that $C$ is the  \emph{discriminant set} of $\vfi$. Note that $\f_C=\emptyset$ if $\# C\not\in 2\bZ$. We set
\[
\f_C(\mu \geq m):=\bigl\{\, (\Gamma,\vfi)\in  \f_C;\;\;\mu(\Gamma)\geq m\,\bigr\},\;\;F(C,m):=\#\f_C(\mu \geq m),
\]
\[
\f_C(\mu = m):=\bigl\{ \,(\Gamma,\vfi)\in  \f_C;\;\;\mu(\Gamma)= m\,\bigr\},\;\;f(C,m):=\# \f_C(\mu = m).
\]
Observe that $F(C,m)=F(C',m)$ if $\#C =\# C'$. We set
\[
\f_n:=\f_{[2n+2]},\;\;\f_n(m):=\f_{[2n+2]}(\mu=m).
\]
 If $\Gamma\in\f_n(m)$ then the level $m+1$ vertex of $\Gamma$ is a node and  thus it can only be of two types:  negative type ($-$) if $m+1$ has two lower neighbors, and positive type ($+$),  if $m+1$ has a unique lower neighbor and thus two higher neighbors. Correspondingly, we obtain a partition
\[
\f_n(m)=\f_n(m)^+\sqcup \f^-_n(m).
\]
We set $F^\pm_n(m):=\#\f^\pm_n(m)$ so that
\begin{equation}
 F_n(m)=F^+_n(m)+F^-_n(m).
\label{eq: pm}
\end{equation}
We discuss two cases.

\smallskip

\begin{figure}[ht]
\centerline{\epsfig{figure=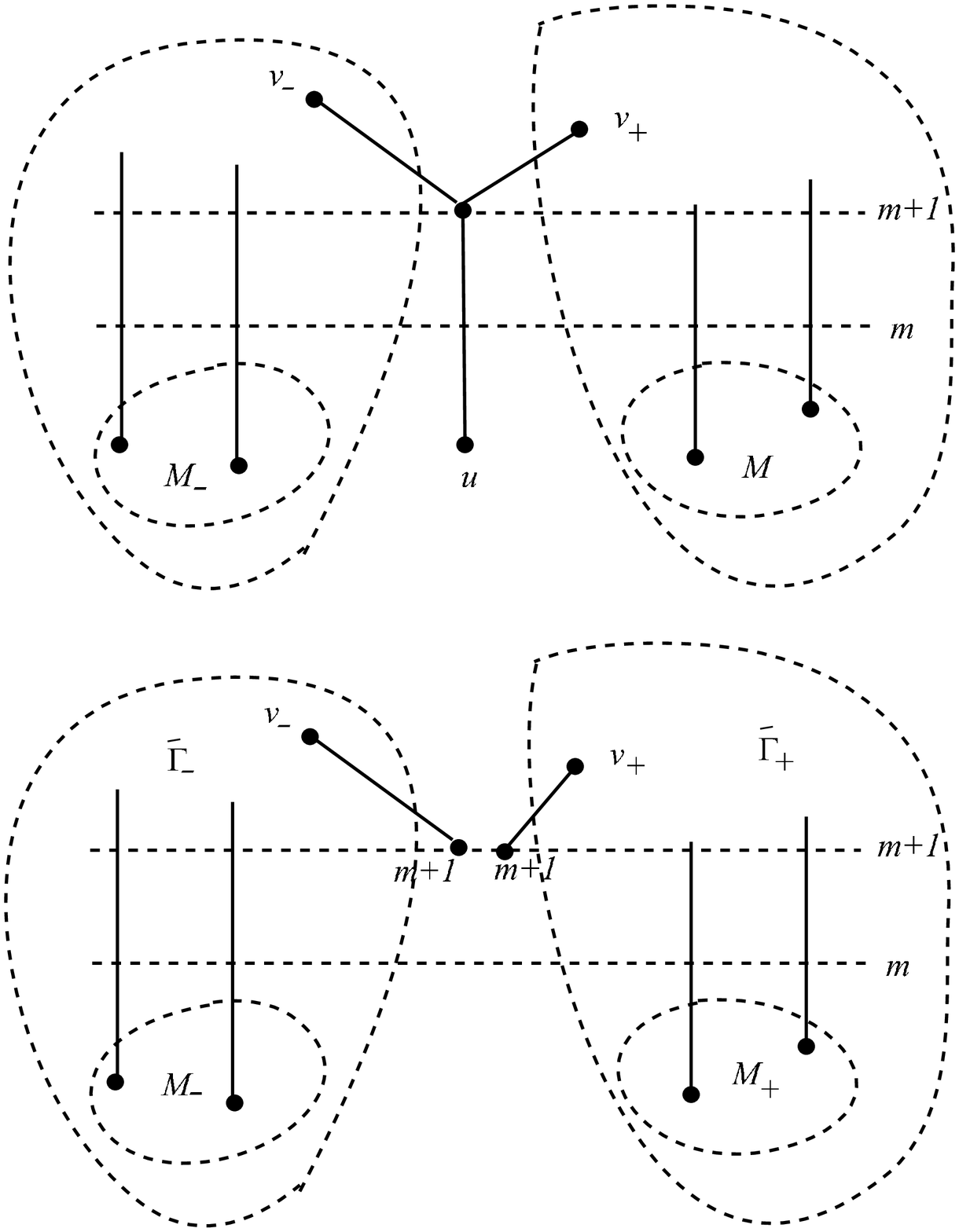,height=3.9in,width=2.9in}}
\caption{\sl Cutting a Morse tree along the lowest $H_1^+$-handle.}
\label{fig: 8}
\end{figure}

\noindent $\mathbf{C}_+$. Suppose $(\Gamma,\ell)\in \f_n(m)^+$. Denote by $v_\pm$ the two higher neighbors of $m+1$  and by $u$ the lower neighbor. We set
\[
k_\pm =\ell(v^\pm),\;\; k=\ell(u),\;\;k^+>k^->m+1>k.
\]
Denote by $\Gamma_{u,m}$ the graph obtained  from $\Gamma$ by removing the vertices $u,m+1$ and the edges  at these points.  Denote by $\Gamma^\pm_{u,m}$ the component containing $v_\pm$ and by $\cv_\pm$ its vertex set. Define $\bar{\Gamma}_\pm$ by setting
\[
\cv(\bar{\Gamma}_\pm)= \cv(\Gamma^\pm_{u,m})\sqcup \{r_\pm\},
\]
where $r_\pm$  has only one neighbor in $\bar{\Gamma}^\pm$, the vertex $v_\pm$. Now define (see Figure \ref{fig: 8})
\[
\vfi_\pm:\cv(\bar{\Gamma}_\pm)\ra \bR,\;\;\ell_\pm(r_\pm)=m+1,\;\;\\vfi_\pm(w)=\ell(w),\;\;\forall w\in\cv(\Gamma^\pm_{u,v}).
\]

Denote by $M_\pm$ the subset of  $\{1,\cdots, m\}$ consisting of vertices which belong to the component $\Gamma_{u,m}^\pm$.   Let $m_\pm$ denote the cardinality of $M_\pm$. We have
\[
M_+\cup M_-= M_k=\{1,\cdots, m\}\setminus \{\ell(u)\} \Longrightarrow m_++m_-=  m-1.
\]
Clearly $M_\pm\subset V_\pm$.  If we set
\[
A_\pm= \cv_\pm \setminus M_\pm
\]
we deduce that $\{A_-,A_+\}$ is a partition of $L_{m+2,n}=\{m+2,\cdots, 2n+2\}$.

$(\bar{\Gamma}_\pm,\vfi_\pm)$  is an unormalized    Morse tree such that its first $m_\pm+1$ vertices are local minima.   The discriminant set of $\vfi_\pm$ is
\[
{C}_\pm=\cv_\pm\cup\{m+1\}= A_\pm\sqcup M_\pm\sqcup\{m+1\},
\]
so that
\[
\bar{\Gamma}_\pm\in\f_{{C}_\pm}(\mu \geq m_\pm +1).
\]
Observe that $(\Gamma, \ell)$ is  determined by the integer $k\in \{1,\cdots, m\}$, the partition
\[
M_k=M_+\sqcup M_-,\;\;\# M_\pm=m_\pm,
\]
the partition
\[
L_{m+2,n}=A_-\sqcup A_+
\]
and the choices of    Morse trees $(\Gamma_\pm,\vfi_\pm)$  with discriminant sets
\[
C_\pm = M_\pm\sqcup A_\pm \sqcup\{m+1\},
\]
such that $\mu(\Gamma_\pm) \geq m_\pm+1$.  The order of $\Gamma_\pm$ is
\begin{equation}
n_\pm =\frac{1}{2}(\# A_\pm + m_\pm -1).
\label{eq: npm}
\end{equation}
If $a_\pm :=\# A_\pm$ we deduce  $a_\pm +m_\pm \equiv 1\mmod 2$ and
\[
a_\pm= 2n_\pm -m_\pm+1,\;\;n_++ n_-=n-1.
\]
We have thus produced a  surjection
\begin{equation}
\vfi:\prod_{k=1}^m\prod_{(C_-^k,C_+^k)}\f_{C_-^k}(\mu\geq m_-+1)\times \f_{C_+^k}(\mu\geq m_++1)\Lra \f_n(m),
\label{eq: surj}
\end{equation}
where the ordered pairs $(C_-^k,C_+^k)$ satisfy
\[
C_\pm^k= M_\pm^k\sqcup A_\pm \sqcup\{m+1\},\;\; M_-^k\sqcup M_+^k=\{1,\cdots, m\}\setminus\{ k\}
\]
\[
A_-\sqcup A_+=\{m+2,\cdots, 2n+2\} ,\;\;\# A_\pm\equiv \# M_\pm^k +1\mmod 2.
\]
We have a \emph{fixed-point-free}  involution on the left-hand-side of (\ref{eq: surj}) defined by the   bijections formally induced by the sign changes $+\longleftrightarrow-$
\[
\f_{C_-^k}(\mu\geq m_-+1)\times \f_{C_+^k}(\mu\geq m_++1)\Lra\f_{C_+^k}(\mu\geq m_++1)\times \f_{C_-^k}(\mu\geq m_-+1).
\]
The fibers of $\vfi$ are precisely the orbits of this involution. We deduce
\[
f^+_n(m)=\frac{1}{2}\sum_k\sum_{(C_-^k,C_+^k)} \#\f_{C_-^k}(\mu\geq m_-+1)\times \#\f_{C_+^k}(\mu\geq m_++1).
\]
The integer $1\leq k\leq m $ can be chosen in $m$ different ways, the partition $\{M_-,M_+\}$ of $\{1,\cdots,m\}\setminus k$ can be chosen in $\binom{m-1}{m_-,m_+}$ ways, then we split $n-1=n_-+n_+$ and  we can choose the partition $\{A_-,A_+\}$ of $\{m+2,\cdots, 2n+2\}$ in $\binom{2n-m+1}{a_-,a_+}$ ways, $a_\pm=2n_\pm-m_\pm+1$,  and finally, the Morse trees  $\Gamma_\pm\in\f_{C_\pm}(m_\pm+1)$ can be chosen in $F_{n_\pm}^+(m_\pm+1)$  ways.  We deduce
\[
F^+_n(m)
\]
\begin{equation}
=\frac{m}{2}\sum_{m_0+m_1=m-1}\;\sum_{n_0+n_1=n-1}\binom{m-1}{m_0} \binom{2n-m+1}{2n_0-m_0+1}F_{n_0}(m_0 +1)F_{n_1}(m_1+1),
\label{eq: f+}
\end{equation}
where for typesetting reasons we  used the notations $m_{0/1}=m_\pm$, $n_{0/1}=n_\pm$.

\begin{figure}[ht]
\centerline{\epsfig{figure=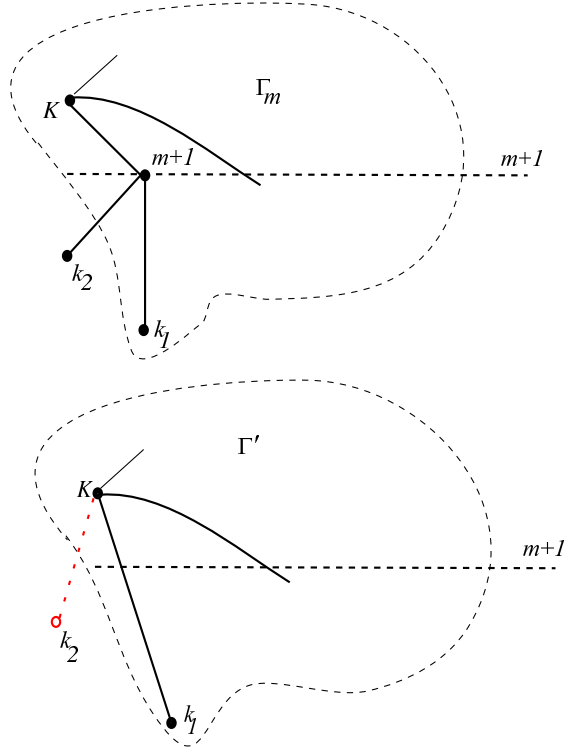,height=3.9in,width=2.9in}}
\caption{\sl Cutting a Morse tree along the lowest $H_1^-$-handle.}
\label{fig: 9}
\end{figure}

\noindent $\mathbf{C}_-$. Suppose $(\Gamma,\ell)\in \f^-_n(m)$. In  this case the vertex $m+1$ has two lower neighbors $1\leq k_1<k_2\leq m$ and a higher neighbor $K>m+1$. Remove the vertex $k_2$ and the unique edge  at $k_2$. The resulting graph $\Gamma_m$ is connected (see Figure \ref{fig: 9}).   We can now produce a Morse tree $(\Gamma', \ell')$ of order $n-1$ as follows (see Figure \ref{fig: 9}). As a graph, $\Gamma'$ is obtained from $\Gamma_m$  by removing the vertex $m+1$ and the two edges $(m+1,K), (m_+,k_1)$ and then connecting $k_1$ and $K$ by an edge. Thus
\[
\cv(\Gamma')=\cv(\Gamma_m)\setminus\{m+1\}\subset \cv(\Gamma)
\]
and we set
\[
\ell':=\ell|_{\cv(\Gamma')}.
\]
The the  discriminant set of $\ell'$ is $[2n+2]\setminus \{m+1,k_2\}$ and $\mu(\Gamma',\ell')\geq m-1$.  The normalized Morse tree $(\Gamma,\ell)$ is uniquely determined by the integer $k_1<k_2\in\{1,2,\cdots, m\}$ and by the un-normalized  Morse tree $(\Gamma',\ell')$ of order $n-1$, with discriminant set $[2n+2]\setminus \{m+1,k_2\}$ and satisfying $\mu(\Gamma')\geq m-1$. We deduce
\begin{equation}
F^-_n(m)= \binom{m}{2}F_{n-1}(m-1).
\label{eq: f-}
\end{equation}
From  (\ref{eq: pm}), (\ref{eq: f+}) and (\ref{eq: f-}) we deduce
\[
f_n(m)= \binom{m}{2}F_{n-1}(m-1)
\]
\[
+ \frac{m}{2}\sum_{m_1+m_2=m-1}\sum_{n_1+n_2=n-1} \binom{m-1}{m_1}\binom{2n-m+1}{2n_1-m_1+1}F_{n_1}(m_1+1) F_{n_2}(m_2+1)
\]
 which is the first equality in  Theorem \ref{th: main}. \qed

\begin{ex} For convenience define $F_n(m)=F_n(1)$ for $m\leq 0$ so that $f_n(m)=0$ for $m\leq 0$. Let $n=0$. Then we have
\[
f_0(1)=1,\;\;f_0(m)=0,\;\;\forall m\neq 1\Longrightarrow F_0(1)=1,\;\;F_0(m)=0,\;\;\forall m>1.
\]
We depicted in Figure \ref{fig: 101} the only Morse function on $S^2$ with  $2=2\times 0+2$ critical points.

\begin{figure}[ht]
\centerline{\epsfig{figure=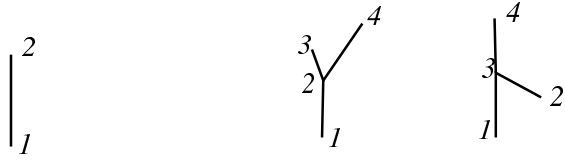,height=0.8in,width=2.9in}}
\caption{\sl Morse functions with $2$ and $4$ critical points.}
\label{fig: 101}
\end{figure}

Let $n=1$. Then Theorem \ref{th: main} predicts
\[
f_1(1)=1,\;\; f_1(2)=F_0(1)=1\Longrightarrow F_1(1)=2,\;\;F_1(2)=1.
\]
The two Morse functions are depicted in Figure \ref{fig: 101}.

For  $n=2$, $m=1$ we have
\[
f_2(1)=\frac{1}{2}\sum_{j=0}^1\binom{4}{2j+1} F_{j}(1)F_{1-j}(1)= 4 F_0(1) F_1(1)=8.
\]
Thus we have $8$ Morse functions with $6$-critical points such that the second critical point is  not a local minimum.  They are depicted in Figure \ref{fig: 10}, where $(S)$:=selfdual.

\begin{figure}[ht]
\centerline{\epsfig{figure=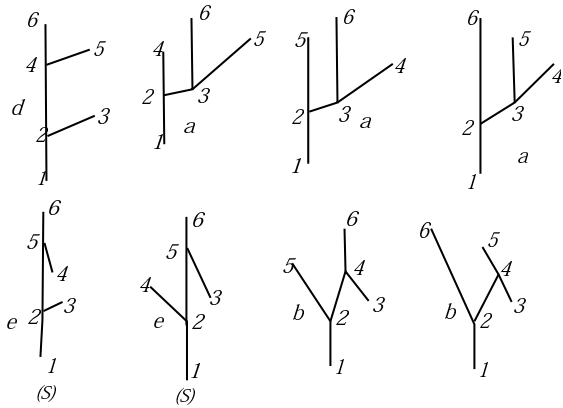,height=2.1in,width=2.9in}}
\caption{\sl Morse functions with  $6$ critical points,  and the second is a saddle point. }
\label{fig: 10}
\end{figure}

For $n=2$, $m=2$ we have
\[
f_2(2)= F_1(1)+\sum_{j=0}^1\binom{3}{j+1} F_{j}(j+1)f_{1-j}(2-j))
\]
\[
= F_1(1)+ \binom{3}{2}F_1(2)F_0(1)+\binom{3}{1} F_0(1)F_1(2)=2+3+ 3=8.
\]
These functions are depicted in Figure \ref{fig: 11}, where $(S)$:=selfdual.
\begin{figure}[ht]
\centerline{\epsfig{figure=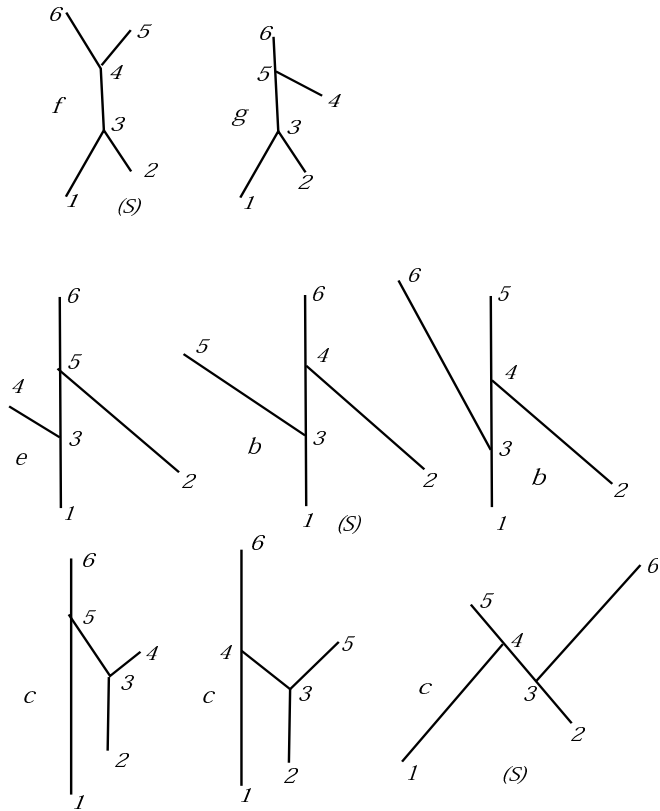,height=4.2in,width=3.4in}}
\caption{\sl Morse functions with $6$ critical points, and the first two are minima. }
\label{fig: 11}
\end{figure}
\begin{figure}[ht]
\centerline{\epsfig{figure=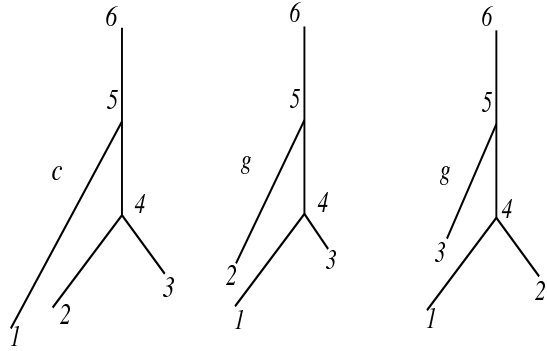,height=1.8in,width=2.8in}}
\caption{\sl Morse functions with $6$ critical points, three minima.}
\label{fig: 12}
\end{figure}

Finally, for $n=2$ and $m=3$ we have $f_2(3)=\binom{3}{2}F_1(2)=3$. These three Morse functions are depicted in   Figure \ref{fig: 12}. Hence
\[
F_2(1)=f_2(1)+f_2(2)+f_2(3)=8+8+3=19,\;\; F_2(2)=11,\;\;F_2(3)=3.
\]
In  Figure \ref{fig: 10}, \ref{fig: 11}, \ref{fig: 12} the small Latin  characters  correspond to the profiles depicted in Figure \ref{fig: 18}.

We list below the  numbers of Morse profiles, homology classes and geometric classes  of Morse functions with $2n+2$  critical points, $ n\leq 9$. The number of geometric equivalence classes were computed using a simple $Maple$ procedure based on the above recurrence. For the reader's convenience we have included this procedure in the Appendix to this paper.
\begin{equation}
\begin{tabular}{||r|r|r|r|r||}\hline
$n$ & $2n+2$& Morse profiles& Homology classes& Geometric classes\\ \hline
0 & 2& 1&1&1\\
1 & 4& 2&2&2\\
2 & 6& 7&10&19\\
3 & 8& 30&70& 428\\
4 & 10& 143&588& 17,746\\
5 & 12& 728& 5,544&1,178,792\\
6 & 14& 3876& 56,628& 114,892,114\\
7 & 16 &21,318&613,470&15,465,685,088\\
8 & 18& 120,175&6,952,660& 2,750,970,320,776\\
9 &  20&690,690& 81,662,152& 625,218,940,868,432\\\hline\hline
\end{tabular}\; .
\label{eq: val}
\end{equation}
\qed
\end{ex}

The  computations in the above example  give an idea of the complexity  of the above recurrence and suggest  that  it can be  better organized.  We do this  in the next section.

\section{Generating functions}
\setcounter{equation}{0}

Consider  the main recurrence formula
\[
f_n(m)= \binom{m}{2}F_{n-1}(m-1)
\]
\[
+\frac{m}{2}\sum_{n_1+n_2=n-1}\sum_{m_1+m_2=m-1}\binom{m-1}{m_1}\binom{2n-m+1}{2n_1-m_1+1}F_{n_1}(m_1+1)F_{n_2}(m_2+1).
\]
Let us first make the change in variables, $m=m+1$. We deduce
\[
f_n(m+1)= \binom{m+1}{2}F_{n-1}(m)
\]
\[
+\frac{m}{2}\sum_{n_1+n_2=n-1}\sum_{m_1+m_2=m}\binom{m}{m_1}\binom{2n-m}{2n_1-m_1+1}F_{n_1}(m_1+1)F_{n_2}(m_2+1).
\]
Next, we introduce new functions
\[
g(m,n)=f_n(m+1),\;\; G(m,n)= F_n(m+1),\;\;0\leq m\leq n.
\]
We deduce $g(m,n)= G(m,n)-G(m+1,n)$. For $n\geq m\geq 1$ we have the equality
\[
G(m,n)-G(m+1,n)=\binom{m+1}{2} G(m-1,n-1)
\]
\[
+\frac{m}{2}\sum_{n_1+n_2=n-1}\sum_{m_1+m_2=m}\binom{m}{m_1}\binom{2n-m}{2n_1-m_1+1}G(m_1,n_1)G(m_2,n_2).
\]
Now  we make the change in variables
\begin{equation}
(m,n)=(x,x+y)\Longleftrightarrow(x,y)=(m,n-m),\;\; H(x,y):= G(m,n).
\label{eq: H}
\end{equation}
With these notations the number of geometric equivalence classes  of Morse functions with $2n+2$ critical points is
\[
F_n(1)=G(0,n)= H(0,n).
\]
Then we have
\[
G(m+1,n)= H(x+1,y-1),\;\;G(m-1,n-1)= H(x-1,y).
\]
For $k=1,2$  we  make the change in variables in the double sum
\[
(m_k,n_k)= (x_k,x_k+y_k)\Longleftrightarrow (x_k,y_k)=(m_k,n_k-m_k).
\]
Then $x_2= m-x_1= x-x_1$,  $y_2= n-1-y_1= x+y-1-y_1$  so that
\[
(m_1,n_1)+(m_2,n_2)= (m,n-1) \Longleftrightarrow (x_1,x_1+y_1)+(x_2,x_2+y_2)= (x,x+y-1)
\]
\[
\Longrightarrow x_1+x_2=x, \;\; y_1+y_2=y-1.
\]
Now observe that in the double sum we need to have
\[
2n_1-m_1+1\leq 2n-m\Longrightarrow 2x_1+y_1+1\leq 2x+y.
\]
\[
2n_2-m_2+1\leq 2n-m\Longleftrightarrow 2(x-x_1)+ (y-1-y_1)+1\leq 2x+y\Longleftrightarrow 0\leq 2x_1 +y_1.
\]
These inequalities are satisfied if and only if
\[
(x_1,y_1)\in R_{x,y-1}:=\bigl\{ (u,v)\in \bZ^2;\;\; 0\leq u\leq x,\;\; 0\leq v\leq y-1\,\bigr\}.
\]
For a point $(x_1,y_1)\in R_{x,y-1}$ we denote by $(\bar{x}_1,\bar{y}_1)$ its  reflection in the center of $R_{x,y-1}$, i.e.
\[
(x_1,y_1)+(\bar{x}_1,\bar{y}_1)=(x,y-1).
\]
The recurrence can now be rewritten as
\begin{multline}
H(x,y)- H(x+1,y-1)=\binom{x+1}{2}H(x-1,y)\\
+\frac{x+1}{2}\sum_{(x_1,y_1)\in R_{x,y-1}}\binom{x}{x_1}\binom{x+2y}{x_1+2y_1+1} H(x_1,y_1)H(\bar{x}_1,\bar{y}_1).
\label{eq: rec}
\end{multline}
We now introduce the  new function
\[
\hat{H}(x,y):=\frac{1}{u(x,y)! v(x,y)!} H(x,y),\;\;u(x,y)=x,\;\;v(x,y)=x+2y+1,
\]
\[
\hat{H}(x,y=-1)=0,\;\;0!:=1.
\]
Observe that
\[
v(x+1,y-1)=x+2y,\;\; v(x-1,y)= x+2y,
\]
\[
v(x_1,y_1)+v(\bar{x}_1,\bar{y}_1)= (x_1+\bar{x}_1)+ 2(y_1+\bar{y_1})+1 =x+2y
\]
We consider two cases.

\noindent {\bf A.} $x>0$. If we  divide both sides of (\ref{eq: rec}) by $x!(x+2y)!$ we deduce that for $x>0$ we have
\[
(x+2y+1)\hat{H}(x,y) - (x+1) \hat{H}(x+1,y-1)
\]
\[
=\frac{x+1}{2} \hat{H}(x-1,y)+\frac{x+1}{2}\sum_{(x_1,y_1)\in R_{x,y-1}} \hat{H}(x_1,y_1)\hat{H}(\bar{x}_1,\bar{y}_1)
\]
{\bf B.} $x=0$. If we divide both sides of (\ref{eq: rec}) by $(2y)!$ we obtain
\[
(2y+1)\hat{H}(0, y)-\hat{H}(1,y-1)=\frac{1}{2}\sum_{y_1=0}^{y-1} \hat{H}(0,y_1)\hat{H}(0,y-1-y_1).
\]
Observe that if  we let $y=0$ in {\bf A} we deduce
\[
\hat{H}(x,0)=\frac{1}{2} \hat{H}(x-1,0)
\]
so that $\hat{H}(x,0)=2^{-x}$. Consider the formal power series
\[
\h(s,t):=\sum_{x,y\geq 0} \hat{H}(x,y) s^x t^{y}.
\]
If we multiply both sides of {\bf A}  and {\bf B} by $s^xt^{y-1}$ and sum over $x\geq 0$, $y\geq 1$ we deduce
\[
\sum_{x\geq 0,y\geq 1} (x+2y+1)\hat{H}(x,y)s^xt^{y-1}- \sum_{x\geq 0,y\geq 1}(x+1) \hat{H}(x+1,y-1)s^xt^{y-1}
\]
\[
= \sum_{x\geq 1,y\geq 1}\frac{x+1}{2} \hat{H}(x-1,y)s^xt^{y-1}+\sum_{x\geq 0, y\geq 1} \frac{x+1}{2}\Biggl(\sum_{R_{x,y-1}}\hat{H}(x_1,y_1)\hat{H}(\bar{x}_1,\bar{y}_1)\Biggr)s^xt^{y-1}.
\]
Make the change in variables $y=y+1$. Then
\[
\sum_{x\geq 0,y\geq 0} (x+2y+3)\hat{H}(x,y+1)s^xt^{y}- \sum_{x\geq 0,y\geq 0}(x+1) \hat{H}(x+1,y)s^xt^{y}
\]
\[
= \sum_{x\geq 1,y\geq 0}\frac{x+1}{2} \hat{H}(x-1,y+1)s^xt^{y}+\sum_{x\geq 0, y\geq 0} \frac{x+1}{2}\Biggl(\sum_{R_{x,y}}\hat{H}(x_1,y_1)\hat{H}(\bar{x}_1,\bar{y}_1)\Biggr)s^xt^y.
\]
Now make the change in variables $x=x+1$ in the third sum.
\[
\sum_{x\geq 0,y\geq 0} (x+2y+3)\hat{H}(x,y+1)s^xt^{y}- \sum_{x\geq 0,y\geq 0}(x+1) \hat{H}(x+1,y)s^xt^{y}
\]
\[
= \sum_{x\geq 0,y\geq 0}\frac{x+2}{2} \hat{H}(x,y+1)s^{x+1}t^{y}+\sum_{x\geq 0, y\geq 0} \frac{x+1}{2}\Biggl(\sum_{R_{x,y}}\hat{H}(x_1,y_1)\hat{H}(\bar{x}_1,\bar{y}_1)\Biggr)s^xt^y.
\]
We obtain
\begin{equation}
\frac{1}{t}\pa_s\bigl(s\h-s\h_{t=0}\,\bigr)+2\pa_t\h-\pa_s\h=\frac{1}{2t}\pa_s\bigl(\,s^2\h-s^2\h_{t=0}\,\bigr)+\frac{1}{2}\pa_s(s\h^2).
\label{eq: temp}
\end{equation}
From the equality
\[
\h_{t=0}=\h(s,0)=\sum_{x\geq 0} 2^{-x}s^x=\frac{2}{2-s}
\]
we obtain
\[
\frac{1}{t}\pa_s(s\h)+2\pa_t\h -\pa_s\h=\frac{1}{2}\pa_s(s\h^2)+\frac{1}{2t}\pa_s(s^2\h)+\frac{1}{t}\pa_s\frac{2s-s^2}{2-s}.
\]
Multiplying both sides by $t$ we obtain
\begin{equation}
\pa_s(s\h)+2t\pa_t\h -t\pa_s\h=\frac{t}{2}\pa_s(s\h^2)+\frac{1}{2}\pa_s(s^2\h)+1,\;\;\h(s,0)=\frac{s}{2-s}.
\label{eq: temp1}
\end{equation}

The above equality is a first order quasilinear p.d.e.  However the   initial condition is  characteristic (see \cite[II.1]{CH}) and thus the above initial value problem   cannot be solved using the method of characteristics. To remove the singularities of this equation we blow it up via   a monoidal  change of coordinates
\[
s=uv,\;\;t=v^2,\;\;\xi(u,v)=v\h(uv,v^2).
\]
Note that
\[
\xi(u,v)=\sum_{x,y\geq 0}H(x,y) \frac{u^xv^{x+2y+1}}{x!(x+2y+1)!}=\sum_{a\geq 0,b\geq 1} H\Bigl(\,a,\frac{b-a-1}{2}\,\Bigr) \frac{u^av^b}{a!\,\cdot\, b!}.
\]
We deduce
\[
v=t^{1/2},\;\;u=st^{-1/2},\;\;\h=v^{-1} \xi,\;\;\;s\h= u\xi,\;\;s^2\h=u^2v\xi,\;\;s\h^2=uv^{-1}\xi^2,
\]
\[
\pa_s=(\pa_su)\pa_u+(\pa_sv)\pa_v= t^{-1/2}\pa_u=v^{-1}\pa_u,
\]
\[
\pa_t= (\pa_tu)\pa_u+(\pa_tv)\pa_v= -\frac{st^{-3/2}}{2}\pa_u+ \frac{t^{1/2}}{2}\pa_v=-\frac{u}{2v^2}\pa_u+\frac{1}{2v}\pa_v=\frac{1}{2v^2}(-u\pa_u+v\pa_v).
\]
The equation (\ref{eq: temp1}) can now be rewritten as
\[
v^{-1}\pa_u(u\xi)+(-u\pa_u+v\pa_v)(v^{-1}\xi)-v\pa_u(v^{-1}\xi)=\frac{v}{2}\pa_u(uv^{-1}\xi^2)+\frac{1}{2v}\pa_u(u^2v\xi)+1.
\]
Multiplying both sides of the above equality by $v$ we deduce
\[
\pa_u(u\xi)+(-uv\pa_u+v^2\pa_v)(v^{-1}\xi)-v\pa_u\xi=\frac{v}{2}\pa_u(u\xi^2)+\frac{v}{2}\pa_u(u^2\xi)+v,
\]
or equivalently
\[
-v\pa_u\xi+v\pa_v\xi=\frac{v}{2}\pa_u(u\xi^2)+\frac{v}{2}\pa_u(u^2\xi)+v.
\]
Dividing by $v$ we obtain
\[
-\pa_u\xi+\pa_v\xi=\frac{1}{2}\pa_u(u\xi^2)+\frac{1}{2}\pa_u(u^2\xi)+1.
\]
This is a first order quasilinear equation  with canonical form
\begin{equation}
-(1+u\xi+\frac{u^2}{2})\pa_u\xi+\pa_v\xi=\frac{1}{2}\xi^2+u\xi+1.
\label{eq: final}
\end{equation}
The characteristic vector field of this equation is (see \cite[\S 7.E]{Ar2} or \cite[II.1]{CH})
\[
V=-(1+u\xi+\frac{u^2}{2})\pa_u+\pa_v+\bigl(\, 1+u\xi+\frac{1}{2}\xi^2\,\bigr)\pa_\xi.
\]
Consider the curve $s\ra \gamma(s)$  described by the initial conditions
\[
u=s,\;\;v=0,\;\;\xi=\xi(s,0)=0\Longrightarrow\frac{d\gamma}{ds}=\pa_u.
\]
Along this curve we have $V(s,0,0)=-(1+\frac{s^2}{2})\pa_u+\pa_v+\pa_\xi$ which  shows that the initial curve is non-characteristic.

The characteristic curves of (\ref{eq: final}) are the solutions of the system of o.d.e.-s
\begin{equation}
\left\{
\begin{array}{rcl}
\frac{du}{dt} &= & -(1+u\xi+\frac{u^2}{2})\\
&&\\
\frac{dv}{dt} &=& 1\\
&&\\
\frac{d\xi}{dt} &=& 1+u\xi+\frac{1}{2}\xi^2
\end{array}
\right..
\label{eq: char}
\end{equation}
The  graph of $\xi$ is filled by the integral curves of (\ref{eq: char}) with initial points on $\gamma$, i.e.
\begin{equation}
u(0)=s,\;\;v(0)=0,\;\;\xi(0)=0.
\label{eq: char0}
\end{equation}
Consider the function
\[
h(u,\xi):= \frac{1}{2}(u^2\xi+u\xi^2)+u+\xi= (u+\xi)\Bigl(\,\frac{u\xi}{2}+1\,\Bigr).
\]
Observe that the plane curve $t\mapsto (u(t),\xi(t))$ is a   solution of the hamiltonian equation
\[
\left\{
\begin{array}{rcl}
\dot{u}&=&-\pa_\xi h\\
\dot{\xi} &=&\pa_uh
\end{array}
\right..
\]
Since the energy is conserved along the trajectories of a hamiltonian system we deduce $h(u,\xi)=const.$ along the trajectories of (\ref{eq: char}). Thus the solutions  of the initial value problem (\ref{eq: char}) $+$ (\ref{eq: char0}) satisfy
\[
\frac{1}{2}(u^2\xi+u\xi^2)+u+\xi=s,\;\; v=t.
\]
We interpret the first equality as a quadratic equation in $\xi$
\[
\frac{u}{2}\xi^2+\bigl(1+\frac{u^2}{2}\,\bigr)\xi+u-s=0
\]
and we solve for $\xi$
\[
\xi=\frac{-(1+\frac{u^2}{2})+\sqrt{(1+\frac{u^2}{2})^2-2u(u-s)}}{u}.
\]
Above, the choice of plus sign in the quadratic formula is dictated by the fact that the Taylor coefficients of $\xi$ are positive. Thus
\[
1+u\xi+\frac{u^2}{2}=\sqrt{\bigl(1+\frac{u^2}{2}\,\bigr)^2-2u(u-s)}=\sqrt{\frac{u^4}{4}-u^2+2su+1}
\]
so that
\[
\frac{du}{dt}=- \sqrt{\frac{u^4}{4}-u^2+2su+1}.
\]
Set
\[
P_s(u):=\frac{u^4}{4}-u^2+2su+1,\;\;\theta=\theta_s(u):=\int_0^u \frac{d\tau}{\sqrt{P_s(\tau)} }.
\]
Then
\[
\theta_s(u)= C_s-t,\;\; C_s= \theta_s(s)=\int_0^s \frac{d\tau}{\sqrt{P_s(\tau)} }.
\]

\begin{remark}We can write $u$ explicitly as a function of $s$ and $t$  using the Weierstrass formula \cite[\S 20.6, Ex.2]{WW}
\[
u= \Phi_s(C_s-t),\;\;\Phi_s=\theta_s^{-1}=\frac{ \wp'(z) +s\wp(z)}{2\bigl(\,\wp(z) +1\,\bigr)^2-\frac{1}{8} },
\]
where $\wp(z)=\wp_s(z)$ is the Weierstrass function with parameters $g_2=\frac{13}{4}$, $g_3= \frac{3}{4}-s^2$.\qed
\end{remark}

 From the equality $0=\theta_s(u=0)=C_s-t$ we deduce that the  curve $u=0$ admits the parametrization
\[
t=C_s=\theta_s(s).
\]
Now observe that
\[
\xi= \frac{-(1+\frac{u^2}{2})^2+(1+\frac{u^2}{2})^2-2u(u-s)}{u\bigl((1+\frac{u^2}{2})+ \sqrt{(1+\frac{u^2}{2})^2-2u(u-s)}\bigr)}=\frac{2s}{((1+\frac{u^2}{2})+ \sqrt{(1+\frac{u^2}{2})^2-2u(u-s)}}.
\]
Hence for $u=0$ we have
\[
\sum_{y\geq 0}\frac{H(0,y)}{(2y+1)!}v^{2y+1}=\xi(0,v)=\xi(0,t)=s,\;\;t=\theta_s(s).
\]
We have thus proved the  following result.

\begin{theorem}  Denote by $\xi_{n}$ the number of geometric equivalence classes of Morse functions on $S^2$ with $2n+2$ critical points and set
\[
\xi(t) = \sum_{n\geq 0}\xi_{n}\frac{t^{2n+1}}{(2n+1)!}.
\]
Then $\xi$ is the compositional inverse of the  function
\[
\theta=\theta(s)=\int_0^s\frac{d\tau}{\sqrt{\frac{\tau^4}{4}-\tau^2+2s\tau+1}},
\]
i.e. $\xi(\theta(s))=s$.\qed
\label{th: main1}
\end{theorem}

\begin{ex} The Taylor coefficients of $\xi(t)$ can be in principle computed  from the above formula via the Lagrange inversion formula although this procedure is not as effective as the recurrence in Theorem \ref{th: main}.  However, we want to test the validity of Theorem \ref{th: main1} on  special cases.

As in \cite{St2}, for every formal power series $f$ in the variable $x$ we denote by $[x^n]f$ the  coefficient of $x^n$ in the expansion of $f$. The Lagrange inversion formula \cite[Thm.5.4.2]{St2} implies
\begin{equation}
[t^5]\xi= \frac{1}{5}[s^4]\Biggl(\frac{s}{\theta(s)}\Biggr)^5.
\label{eq: lagr}
\end{equation}
We write
\[
P_s= 1+r_s(t),\;\; r_s(t)=2st-t^2+\frac{t^4}{4}.
\]
Then we have a binomial  expansion
\[
P_s(t)^{-1/2}= 1-\frac{1}{2}r_s(t)+\frac{1\cdot 3}{2^2\cdot 2!}r_s(t)^2-\frac{1\cdot 3\cdot 5}{2^3\cdot 3!} r_s(t)^3+\frac{1\cdot 3\cdot 5\cdot 7}{2^4\cdot 4!}r_s(t)^4+\cdots
\]
Integrating this equality   with respect to $t\in[0,s]$ we deduce
\[
\theta(s)=\int_0^sP_s(t)^{-1/2} dt =s\Bigl(\, 1- \underbrace{(\frac{1}{3}s^2-\frac{7}{40}s^4+\frac{3}{28}s^6+\cdots)}_{=:q(s)} \;\;\Bigr).
\]
Then
\[
\frac{s}{\theta(s)}=\frac{1}{1-q(s)},\;\;\Biggl(\frac{s}{\theta(s)}\Biggr)^5= 1+5 q(s) +5\cdot 6\frac{q(s)^2}{2!}+\cdots.
\]
We deduce that
\[
[s^4]\Biggl(\frac{s}{\theta(s)}\Biggr)^5= 5[s^4]q(s)+15[s^4]q(s)^2 =-\frac{35}{40}+\frac{15}{9}=\frac{19}{24}\Longrightarrow 5![t^5]\xi=19.
\]
This agrees with the value computed in (\ref{eq: val}).

Similarly, we have
\[
[t^7]\xi= \frac{1}{7}[s^6]\Biggl(\frac{s}{\theta(s)}\Biggr)^7
\]
and we deduce
\[
[s^6]\Biggl(\frac{s}{\theta(s)}\Biggr)^7= 7[s^6]q(s)+ 7\cdot 8[s^6] \frac{q(s)^2}{2!}+7\cdot 8\cdot 9[s^6]\frac{q(s)^3}{3!}=\frac{107}{180}\Longrightarrow 7![t^7]\xi= 428.
\]
This too agrees with the value  found in (\ref{eq: val}). \qed
\end{ex}

\section{On the topological equivalence problem}
\setcounter{equation}{0}

We were not able to find a computationally satisfactory recurrence for the number  of topological equivalence classes of Morse functions but we could still  describe some interesting   combinatorial structures on this set.

A regular sublevel set of a Morse function on a sphere is a disjoint union of $2$-spheres with  small open disks removed.

The topology of such a disjoint  union of holed spheres is encoded by a partition $\pi$,  i.e.    a  decreasing function
\[
\pi:\bZ_{> 0}\ra \bZ_{\geq 0},\;\; \pi(i) \geq \pi(i+1),\;\;\forall i> 0,
\]
such that $\pi(i)=0$ for all $i\gg 0$.  The length of the partition  is the nonnegative integer.
\[
\ell(\pi) =\max \{ i;\;\;\pi(i)>0\} -1.
\]
 The  weight of the partition is the integer
 \[
 |\pi|= \sum_{i>0}\pi(i).
 \]
 If $n=|\pi|$ we say that $\pi$ is a partition on $n$. We denote by $\p$ the set of all partitions of nonnegative weight and by $\p_n$  the set of partitions of  weight $n$. ${\bf O}$ denotes the unique  partition of weight $0$.

  To a partition $\pi$ of weight $w$ and length $\ell$ there corresponds a disjoint union  of  $\ell$  holed sphere,  $\pi(1)$ holes on the first sphere, $\pi(2)$ holes on the second sphere etc.   Observe that the weight of the partition is equal to the number of  boundary components of the sublevel set.

 We will describe   partitions by finite  sequences of positive integers, where two sequences are to be considered equivalent if one can be obtained  from the other by a permutation. For example $\{4,4,1\}=\{4,1,4\}$ denotes a  union of two spheres with  $4$-holes and a sphere with $1$ hole.  Following a longstanding tradition we will also use exponential notation to indicate partitions. Thus
 \[
 h_1^{m_1}\cdots h_k^{m_k}=\{\,\underbrace{h_1,\cdots,h_1}_{m_1},\cdots,\underbrace{h_k,\cdots,h_k}_{m_k}\,\}.
 \]
The possible transitions from one sublevel set to the next correspond to four simple operations on partitions. Given a partition
\[
h_1\geq h_2\cdots \geq h_\ell >0
\]
these operations are
\[
\{h_1,\cdots, h_\ell\}\stackrel{H_0}{\longmapsto}\{h_1,\cdots, h_\ell, 1\},
\]
\[
\{h_1,\cdots,h_i,\cdots, h_\ell\}\stackrel{H_1^+}{\longmapsto}   \{h_1,\cdots,h_i+1,\cdots, h_\ell\},
\]
\[
\{h_1,\cdots,h_i,\cdots,h_j,\cdots  h_\ell\}\stackrel{H_1^+}{\longmapsto} \{h_1,\cdots,h_{i-1}, h_i+h_j-1,h_{i+1},\cdots,h_{j-1},h_{j+1}\cdots  h_\ell\},
\]
\[
\{h_1,\cdots,h_i,\cdots, h_\ell\}\stackrel{H_2}{\longmapsto}   \{h_1,\cdots,h_i-1,\cdots, h_\ell\},\;\;(h_i>1).
\]
We will refer  to these  transitions as \emph{moves}. Note that the $H_0$ move increases the length and the weight by $1$, the $H_1^+$ move increases the weight by $1$ but preserves the length, the $H_1^-$ move decreases the weight and the length by $1$, while the $H_2$ move  decreases the weight by one but preserves the length.

We introduce a simplified set of moves. Consider a partition
\[
h_1\geq h_2 \geq \cdots\geq h_\ell \geq 0.
\]
The simplified collection of moves is

\smallskip

\noindent $\bullet$ Type $U$ (up) move
\[
\{h_1,\cdots,h_i,\cdots, h_\ell\}\stackrel{U}{\longmapsto} \{h_1,\cdots,h_i+1,\cdots, h_\ell\},
\]
\noindent $\bullet$ Type $D$ (down) move
\[
\{h_1,\cdots,h_i,\cdots, h_\ell\}\stackrel{D}{\longmapsto}   \{h_1,\cdots,h_i-1,\cdots, h_\ell\},\;\;h_i>0.
\]
\noindent $\bullet$ Type $X$-move\footnote{I call this the $X$ move since at this moment I don't know what to make of it.}
\[
\{h_1,\cdots,h_i,\cdots,h_j, \cdots  h_\ell\} \stackrel{X}{\longmapsto} \{h_1,\cdots,h_i+h_j-1,\cdots,h_{j-1} ,h_{j+1}, \cdots  h_\ell\},\;\; h_i\geq h_j>1.
\]

We can use these moves to produce a  directed graph with vertex set $\p$.   The procedure is very simple. For every  move
\[
\p\ni \pi\ra \pi'\in\p
\]
draw an arrow directed from $\pi$ to $\pi'$.  We denote by $\hat{\p}$ this graph.  Then the number of  topological equivalence classes  of Morse functions with $2n+2$ critical points  is equal to the number of  directed paths of length $2n$ from the partition $1$ to itself. Equivalently, its is the number of  directed paths of length $2n+2$ from ${\bf O}$ to itself. We denote this number by $T_{2n+2}$.

We say that a path in $\p$ is \emph{constructive}  if it consists
only of the $U$-moves.  A path is called \emph{simple} if it consists only of the $U$, $D$ moves. Given $\pi,\pi'\in\p$ we write $\pi\prec\pi'$ if there exists a constructive path from $\pi$ to $\pi'$. We  have the following elementary fact whose proof is left to the reader.

\begin{lemma}   The following statements are equivalent.

\noindent (a) $\pi\prec\pi'$.

\noindent (b) $\pi\subseteq \pi'$  i.e. $|\pi|<|\pi'|$ and $\pi(i)\leq \pi'(i)$, $\forall i>0$.

\end{lemma}

The  order relation  $\pi\subseteq  \pi'$  on $\p$ is known as the Young ordering (see \cite[Chap.7]{St2}). For every $\mu\subseteq \lambda$  we denote by $Z_n(\lambda,\mu)$ the number of $r$-chains   from $\mu$ to $\lambda$, i.e. the number  sequences $(\mu_1,\cdots,\mu_r)$ satisfying
 \[
 \mu\subsetneq \mu_1\subsetneq\cdots \subsetneq\mu_r\subsetneq\lambda.
 \]
 we have the  formula of Kreweras \cite{Kr} (see  also \cite[\S 2.4]{Be} or \cite[Exer. 3.63]{St1}))
\[
Z_{r}(\mu,\lambda)=\det\left[\binom{\lambda_i-\mu_j+r}{i-j+r}\right]_{1\leq i, j\leq \ell},\;\;\ell=\ell(\lambda).
\]
Then if $r=|\lambda|-|\mu|-1$ then $Z_r(\mu,\lambda)$ is precisely the number of constructive paths from $\mu$ to $\lambda$. We  denote this number by $C(\mu,\lambda)$. When $\mu={\bf O}$ this  formula simplifies considerably. More precisely  $C({\bf O}, \lambda)$ is the number  of \emph{standard  Young tableaux} (SYT) of shape $\lambda$ (see \cite[Prop. 7.10.3]{St2}).  More precisely we have the \emph{hook-length formula}
\[
C(\lambda)=C({\bf O},\lambda)= \frac{|\lambda|!}{\prod_{u\in \lambda}h(u)},
\]
where the above product  is taken over all the cells $u$ of the Young diagram of $\lambda$ and $h(u)$ is the hook length of $u$,
\[
h(u) =1+\mbox{ \# \{cells to the right of $u$\} +\# \{cells below $u$\}}.
\]

\begin{figure}[ht]
\centerline{\epsfig{figure=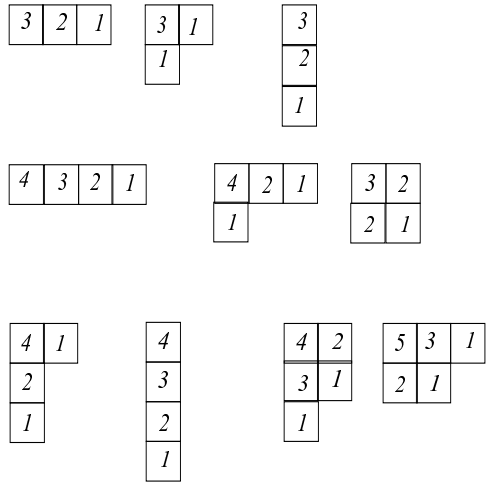,height=2.5in,width=2.5in}}
\caption{\sl Young diagrams of weight $3$ and $4$ and $5$. }
\label{fig: 21}
\end{figure}

\begin{ex}  In Figure \ref{fig: 21} we have  depicted all the Young diagrams of weight $3$ and $4$ and two Young  diagrams of weight $5$. The numbers inside the cells  the hook lengths. For a  Young diagram $\lambda$ we set
\[
h(\lambda)=\prod_{u\in \lambda}h(u)
\]
and we will refer to it as the \emph{hook weight}. In the table below we have recorded the hook weight of the diagrams depicted in Figure \ref{fig: 21}.
\[
\begin{tabular}{||l|r||}\hline
$\lambda$ & $h(\lambda)$ \\ \hline
$3$ & 6\\
$2,1$ &3 \\
$1^3$ &6\\\hline
$4$ &$24$\\
$3,1$ &$8$\\
$2^2$ &$12$\\
$2,1^2$ & $8$\\
$1^4$ & $24$ \\ \hline
$2^2,1$ & $24$\\
$3,2$ & $30$\\ \hline\hline
\end{tabular}
\]
\qed
\end{ex}

Let us point out that  the number of simple paths of length $2n+2$ from ${\bf O}$ to itself has been computed in \cite[ Eq. (39)]{St3} and it is
\[
S_{2n+2}= \frac{(2n+2)!}{2^{n+1}(k+1)!}=1\cdot 3\cdots(2n+1).
\]
In particular we deduce
\[
T_{2n+2}\geq 1\cdot 3\cdots(2n+1).
\]
\begin{ex} (a) Since the $X$-moves   do not affect   partitions which have only one part $>1$     we deduce that for $n=2$ we have
\[
T_6=S_6= 15
\]
so that there are exactly $15$ topological equivalence classes of Morse functions on $S^2$ with $6$ critical points.

\noindent (b) Note that many  of the edges  of the graph $\hat{\p}$   are edges of the   Hasse diagram of the Young lattice $(\p,\subseteq)$.  Remove  these edges to obtain  a new graph $\p_X$.  There is an arrow $\mu\ra \lambda$ in this graph if and only  $\lambda$ can be reached from $\mu$  by one $X$ move,  but cannot be reached  from $\mu$ by a $D$-move.  For example, $(21)$ can be reached from $(22)$ by a $D$-move  but cannot be reached  by an $X$-move  so they are not connected  in $\p_X$.

\begin{figure}[ht]
\centerline{\epsfig{figure=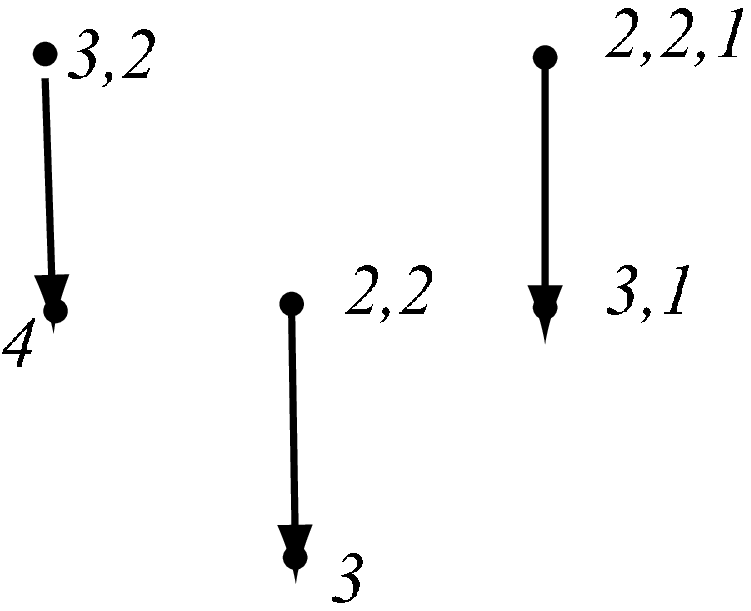,height=1.2in,width=1.5in}}
\caption{\sl The graph $\p_X$. }
\label{fig: 22}
\end{figure}

 Let us compute $T_8$.  A path in $\hat{\p}$ of length  $8$ from ${\bf O}$ back to ${\bf O}$ goes only through partitions of weight $\leq 4$.  If $\gamma$ is  such a path  which    contains an $X$-edge then that edge must be the edge $(2,2)\ra 3$ depicted in Figure $\ref{fig: 22}$.  We  deduce that
\[
T_8= S_8 + C(2,2)\times C(3)= 15\times 7 + 2\times 1= 107.
\]
\noindent (c) Denote by $S_k(\mu| \lambda)$ the number of  simple paths in $\hat{\p}$ of length $k$ from $\mu$ to $\lambda$ and by $T_k(\mu|\lambda)$ the number of all path in $\hat{\p}$  from $\mu$ to $\lambda$.  Upon inspecting Figure \ref{fig: 22} we deduce
\[
T_{10}=S_{10}+ C({\bf O}|2,2,1)\times \bigl( \,C({\bf O}|3,1)+ C({\bf O}|3)\,\bigr)
\]
\[
+ C({\bf O}|3,2)\times \bigl( C({\bf O}|4)+ C({\bf O}|3)\,\bigr)
\]
\[
+ C({\bf O}|2,2)\times \bigl(\, C({\bf O}|3,1)+ C({\bf O}|4)\,\bigr)
\]
\[
= 9\times 105 + 5\times(3+1) + 4\times(1+1)+ 2\times(3+1)=945+20+8+8= 981.\proofend
\]
\end{ex}

\newpage

\appendix
\section{$Maple$ implementation}

We include below  a simple but  far from optimal   $MAPLE$ procedure for computing the number  of Morse functions  based on the recurrence satisfied by the numbers $\hat{H}(x,y)$ described in Section 8.

\begin{verbatim}
>Morse:=proc(a::nonnegint, b:: nonnegint)
>local i,j,k, m, x, y, A;A[0,0]:=1; m:=a+b;
>for k from 1 to m  do
>for y from 0 to b do
>x:= k-y;
>if y=0 then A[x,y]:=1/(2^x) elif x>0 then
>A[x,y]:=(1/(x+2*y+1))*( (x+1)*A[x+1,y-1]+ (1/2)*(x+1)*A[x-1,y]+(1/2)*(x+1)
*add( add(A[i,j]*A[x-i,y-1-j],j=0..y-1),i=0..x)) else
> A[x,y]:=(1/(2*y+1))*( (x+1)*A[x+1,y-1]+ (1/2)*(x+1)
*add( add(A[i,j]*A[x-i,y-1-j],j=0..y-1),i=0..x))
> end if;
> end do;
> end do;
> (a!)*(2*b+1)!*A[a,b];
> end proc:
\end{verbatim}

To compute the number $H(x,y)$  in (\ref{eq: H}) use the command

\begin{verbatim}
>Morse(x,y);
\end{verbatim}

The number of geometric equivalence classes of Morse functions with $2n+2$ vertices is obtained using the command

\begin{verbatim}
>Morse(0,n);
\end{verbatim}

To deal with large $n$ ($n>25$)  modify the  command line

\begin{verbatim}
> (a!)*(2*b+1)!*A[a,b];
\end{verbatim}

to

\begin{verbatim}
>A[a,b];
\end{verbatim}

The procedure will then  generate the number $\frac{1}{(2n+1)!} H(0,n)$.

\end{document}